\documentclass[a4paper,12pt]{amsart}
\usepackage{ae,amsfonts,euscript,enumerate}
\usepackage{amsmath,euscript}
\usepackage{amssymb}
\usepackage{amsthm}
\usepackage{enumerate}
\usepackage{amsfonts}
\usepackage{comment}
\usepackage{colortbl}
\usepackage{a4wide}
\usepackage{amssymb,amsmath,array}

\newtheorem{thm}{Theorem}[section]
\newtheorem{lem}[thm]{Lemma}

\newtheorem{prop}[thm]{Proposition}
\newtheorem{cor}[thm]{Corollary}

\theoremstyle{defn}
\newtheorem{defn}{Definition}

\newtheorem{rem}{Remark}
\newtheorem{quest}{Question}
\newtheorem*{ex}{Example}

\newtheorem*{notn}{Notation}

\newcommand{\PMC}{\mathcal{PM}(C)}

\newcommand{\gray}{\cellcolor[gray]{0.1} }

\newcommand{\gc}{ [ \hspace{-0.65mm} [}
\newcommand{\dc}{]  \hspace{-0.65mm} ]}

\newcommand{\ia}{i,\alpha}

\newcommand{\hc}{\mathcal{H}}

\newcommand{\spec}{{\rm Spec}}

\def\ch{{\mathcal H}}
\def\co{{\mathcal O}}

\def\oqmm13{\co_q(M_{1,3})}
\def\oqm23{\co_q(M_{2,3})}
\def\oqmmn{\co_q(M_{m,n})}

\def\ia{i,\alpha}

\input{xy}
\xyoption{all}

\def\eqref#1{(\ref{#1})}
\def\qed{~\vrule height8pt width 5pt depth -1pt\medskip}

\usepackage[dvips]{graphicx}
\usepackage{amssymb}
\usepackage{amsmath}
\usepackage[latin1]{inputenc}
\newif\ifpdf
       \ifx\pdfoutput\undefined
       \pdffalse 
       \else
       \pdfoutput=1 
       \pdftrue
\fi

\ifpdf
       \else
       \usepackage{graphicx,color,rotating}
\fi
\usepackage{fullpage}
\title[]{An automaton-theoretic approach to the representation theory of quantum algebras}
\author{J.~Bell, S.~Launois, and J.~Lutley}
\thanks{The first author thanks NSERC for its generous support.  The second author research was supported by a Marie Curie European Reintegration Grant within the $7^{\mbox{th}}$ European Community Framework Programme}

\address{Jason Bell\\
Department of Mathematics\\
Simon Fraser University\\
Burnaby, BC V5A 1S6, Canada
}

\email{jpb@math.sfu.ca}

\address{
S. Launois \\
Institute of Mathematics, Statistics \& Actuarial science \\
University of Kent\\
Canterbury, Kent CT2 7NF\\
United Kingdom}

\email{s.launois@kent.ac.uk}

\address{Jamie Lutley\\
Department of Mathematics\\
Simon Fraser University\\
Burnaby, BC V5A 1S6, Canada
}

\email{jlutley@sfu.ca}
\begin{document}
\bibliographystyle{plain}

\begin{abstract}
We develop a new approach to the representation theory of quantum algebras supporting a torus action via methods from the theory of finite-state automata and algebraic combinatorics.  We show that for a fixed number $m$, the torus-invariant primitive ideals in $m\times n$ quantum matrices can be seen as a regular language in a natural way.  Using this description and a semigroup approach to the set of Cauchon diagrams, a combinatorial object that paramaterizes the primes that are torus-invariant, we show that for $m$ fixed, the number of torus-invariant primitive ideals in $m\times n$ quantum matrices satisfies a linear recurrence in $n$ over the rational numbers.  In the $3\times n$ case we give a concrete description of the torus-invariant primitive ideals and use this description to give an explicit formula for the number $P(3,n)$.  
\end{abstract}
\maketitle

\section{Introduction}

For a given infinite dimensional algebra, it is often a very
difficult problem to classify its irreducible
representations. Because of this, Dixmier proposed classifying the
primitive ideals of this algebra as an intermediate step; once one has classified the primitive ideals, one can then for each primitive ideal $P$
try to find an irreducible representation whose annihilator is $P$. 

In this paper, we study primitive ideals in the quantum
world, and in particular primitive ideals of the algebra $\oqmmn$ of generic
quantum matrices. This algebra is a noncommutative deformation of the coordinate ring of the variety of matrices, and our aim is to understand the (noncommutative) geometry of the ``variety of quantum matrices.'' In this spirit, the knowledge of the primitive spectrum of the algebra of generic quantum matrices is of crucial importance;  by analogy with classical algebraic geometry, we think of the primitive ideals of $\oqmmn$ as corresponding to the points of the ``variety of quantum matrices.'' 
 
Several results have already been obtained regarding the primitive ideals of $\oqmmn$. In particular, it is known from work
of Goodearl and Letzter \cite{GL} that the Dixmier-Moeglin equivalence holds for these
algebras. Even better, this algebra supports a ``nice'' torus action and, it follows from the work of Goodearl and Letzter that the number of prime ideals of $\oqmmn$ invariant under the action of this torus $\hc$ is finite. Because of this, the prime spectrum of $\oqmmn$ admits a stratification into finitely many $\hc$-strata. Each $\hc$-stratum is defined by a unique $\hc$-invariant prime ideal---that is minimal in its $\hc$-stratum---and is homeomorphic to the scheme of irreducible subvarieties of a torus. Moreover the primitive ideals correspond to those primes that are maximal in their $\hc$-strata. Hence, it is very important to recognize those $\hc$-invariant prime ideals that are primitive; they correspond to those $\hc$-strata that are homeomorphic to the scheme of irreducible subvarieties of the base field. In other words, the number of $\hc$-invariant prime ideals that are primitive 
gives the number of points that are invariant under the induced action of $\hc$ in the ``variety of quantum matrices.''

In an earlier paper \cite{BLN}, a strategy was developed to recognize the primitive $\hc$-invariant prime ideals.  The basic idea behind this strategy was that for each $\hc$-invariant prime in $\oqmmn$, one can associate an $m\times n$ grid called a \emph{Cauchon diagram} \cite{Cauchon} (or Le-diagrams \cite{Post}) in which squares are colored either black or white with certain restrictions that are described in \S \ref{diagram}.  To each diagram, we then associate a skew-adjacency matrix in a natural way.  Primitivity is equivalent to this matrix being invertible.  In this paper, we extend this strategy by showing that if the determinant is nonzero, it must in fact be a power of 4.  This is very useful, because it means invertibility can be deduced by looking at the determinant mod 3.  Using this fact, and a new approach of the algebra of quantum matrices via the theory of automata, we are able to prove the following theorem.

\begin{thm} Let $m$ be a natural number and let $P(m,n)$ denote the number of $m\times n$ primitive $\hc$-primes in $\mathcal{O}_{q}(M_{m,n})$. Then the sequence $\left\{P(m,n)\right\}^{\infty}_{n=1}$ satisfies a linear recurrence over $\mathbb{Q}$; that is, there exists a natural number $d$ and rational constants $c_{1},c_{2},\ldots,c_{d}$ such that $P(m,n) = \sum_{i=1}^d c_i P(m,n-i)$ for all $n\ge d$.
\label{thm: main1}
\end{thm}

This theorem represents an important step in the understanding of the primitive $\hc$-primes, and a first progress towards a resolution of a conjecture which says that not only should a linear recurrence exist, but that it should have a very specific form \cite[Conjecture 2.9]{BLN}.

In the case of $m = 3$ we are able to give an explicit recurrence and thus prove a conjecture \cite[Conjecture 2.8]{BLN}. Again this is achieved thanks to an automaton theoritic approach.

\begin{thm} Let $n\ge 1$ be a natural number. Then the number of $3\times n$ primitive $\hc$-primes in $\mathcal{O}_{q}(M_{m,n})$ is given by  $$\frac{1}{8} \cdot  \Big(  15\cdot 4^n - 18 \cdot 3^n +13 \cdot 2^n - 6\cdot(-1)^n + 3\cdot (-2)^n \Big).$$
\label{thm: main2}
\end{thm} 

In order to obtain these results we introduce a semigroup structure on the set of diagrams---a diagram is just an $m\times n$ grid with each square colored black or white---and show that the Cauchon diagrams form a regular subset (in an automaton-theoretic sense).  Cauchon/Le-diagrams have been extensively studied recently (see for instance \cite{Cor4,Jo,Ste}) because of their links with other areas: total positivity (see \cite{Post,Wi}), quantum algebras (see \cite{Cauchon}), and the partially asymmetric exclusion process (see \cite{Cor1,Cor2,Cor3}). The results obtained here on Cauchon diagrams are of intrinsic interests. \\

The paper is organized as follows. In \S \ref{diagram}, we give background on Cauchon diagrams and their associated skew-adjacency matrices in \S \ref{Pfaffian}, we recall the basic facts about Pfaffians, which we use to compute the determinants of skew-adjacency matrices; we also introduce the notion of a decomposition of a diagram, and show how this simplifies many conjectures.  In \S \ref{semigroup}, we show that the diagrams can be endowed with a semigroup structure and determine the relation between this semigroup and primitivity. In \S \ref{excess}, we use representation theory to give isomorphisms, which we will later use in proving Theorem \ref{thm: main2}.  In \S \ref{automaton} we give some background in finite state automata and use this to prove Theorem \ref{thm: main1}.  Finally, in \S \ref{enumeration}, we prove Theorem \ref{thm: main2}.

Throughout this paper, we use the following notation and conventions. Given a statement $S$ that is either true or false, we let 
\begin{equation}
\delta(S) \ = \ \left\{ \begin{array}{ll}
 1 & {\rm if ~}
 S
 {\rm ~ is~ true;} \\
0 & {\rm otherwise}. \end{array} \right. 
\end{equation}
We write $\mathbb{Z}_m$ to denote the ring $\mathbb{Z}/m\mathbb{Z}$.
We note that since $(-1)^k$ is constant on cosets of $2\mathbb{Z}$ in $\mathbb{Z}$, we define $(-1)^{a+2\mathbb{Z}}$ to be $(-1)^a$.  Similarly, for the sake of convenience we will sometimes write $a\equiv b~(\bmod\, m)$ when one of $a$ or $b$ is an element of $\mathbb{Z}_m$; by this we just mean that any element in the equivalence class of $a$ is congruent to $b$ mod $m$.  We use bold font to represent vectors and ${\bf 0}$ and ${\bf 1}$ respectively denote the zero vector and the vector whose components are all $1$. 

\section{Diagrams}\label{diagram}
In this section we give the basic background on diagrams and matchings, which are the key component of the combinatorial criterion for an $\hc$-invariant prime ideal in $\oqmmn$ to be primitive.  We begin with a definition.
\begin{defn} {\em
We call an $m\times n$ grid consisting of $mn$ squares in which each square is colored either black or white an $m\times n$ \emph{diagram}. We say that an $m\times n$ diagram with $k$ white squares is \emph{labelled} with labels $\left\{i_{1},i_{2},\ldots,i_{k}\right\}$ if each white square is assigned a number subject to the rules:

\begin{enumerate}
\item The labels are strictly increasing along rows;
\item If a white square is in a higher row than another then its label is the smaller of the two.
\end{enumerate} }
\end{defn}
\begin{figure}[h]
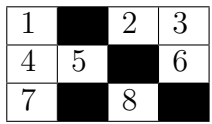

\label{fig: 13}
\vskip 2mm
\begin{tabular}{|p{0.30cm}|p{0.30cm}|p{0.30cm}|p{0.30cm}|}

\hline
$1$ & \gray & $2$ & $3$ 
 \\
\hline
 $4$ & $5$ & \gray & $6$  
 \\
\hline
$ 7$ &\gray & $8$  & \gray\\
\hline
\end{tabular}
\vskip 2mm 
\label{fig: diagram}
\caption{An example of a labelled $3\times 4$ diagram with labels $\left\{1,2,3,4,5,6,7,8\right\}.$}
\end{figure}

We are interested in a special subcollection of diagrams which are known as \emph{Cauchon diagrams}.

\begin{defn} {\em
We say that an $m\times n$ diagram is a \emph{Cauchon diagram} if whenever a square is colored black, either every square in the same row that is to its left is also black or every square in the same column that is above it is also black. We say that a Cauchon diagram is \textit{labelled} if it is a Cauchon diagram and a labelled diagram.
} \end{defn}
\begin{figure}[h]
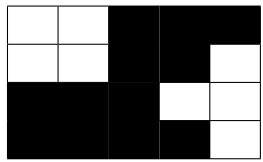

\label{fig: 2}
\vskip 2mm
\begin{tabular}{|p{0.30cm}|p{0.30cm}|p{0.30cm}|p{0.30cm}|p{0.30cm}|}
\hline
& & \gray & \gray &\gray 
 \\
\hline
& &\gray & \gray & 
 \\
\hline
\gray& \gray &\gray &   & 
\\
\hline
\gray& \gray &\gray & \gray  & 
\\
\hline
\end{tabular}
\vskip 2mm
\caption{A $4\times 5$ Cauchon diagram}
\end{figure}
Cauchon diagrams are named for their discoverer, who showed that there is a natural bijection between the collection of $m\times n$ Cauchon diagrams and the $\hc$-primes in $\mathcal{O}_{q}(M_{m,n})$, the ring of $m\times n$ quantum matrices \cite{Cauchon}. Although we won't work directly with the algebra of quantum matrices and its $\hc$-primes, we briefly recall their defintions.

Let $\mathbb{K}$ be a field and $q$ be a nonzero element of $\mathbb{K}$ that is not a root of unity. The algebra $\mathcal{O}_{q}(M_{m,n})$ of quantum matrices is the algebra it is
the $\mathbb{K}$-algebra generated by the $m \times n $ indeterminates
$Y_{\ia}$, $1 \leq i \leq m$ and $ 1 \leq \alpha \leq n$, subject to the
following relations:\\ \[
\begin{array}{ll}
Y_{i, \beta}Y_{i, \alpha}=q^{-1} Y_{i, \alpha}Y_{i ,\beta},
& (\alpha < \beta); \\
Y_{j, \alpha}Y_{i, \alpha}=q^{-1}Y_{i, \alpha}Y_{j, \alpha},
& (i<j); \\
Y_{j,\beta}Y_{i, \alpha}=Y_{i, \alpha}Y_{j,\beta},
& (i <j,  \alpha > \beta); \\
Y_{j,\beta}Y_{i, \alpha}=Y_{i, \alpha} Y_{j,\beta}-(q-q^{-1})Y_{i,\beta}Y_{j,\alpha},
& (i<j,  \alpha <\beta). 
\end{array}
\]

It is well known that $\mathcal{O}_{q}(M_{m,n})$ can be presented as an iterated Ore extension over
$\mathbb{K}$, with the generators $Y_{\ia}$ adjoined in lexicographic order.
Thus the ring $\mathcal{O}_{q}(M_{m,n})$ is a noetherian domain. Moreover, since $q$ is not a root of unity, it follows from
\cite[Theorem 3.2]{gletpams} that all prime ideals of $R$ are completely
prime.

It is easy to check
that the group $\hc:=\left( \mathbb{K}^* \right)^{m+n}$ acts on $\mathcal{O}_{q}(M_{m,n})$ by
$\mathbb{K}$-algebra automorphisms via:
$$(a_1,\dots,a_m,b_1,\dots,b_n).Y_{\ia} = a_i b_\alpha Y_{\ia} \quad {\rm
for~all} \quad \: (\ia)\in \gc 1,m \dc \times \gc 1,n \dc.$$ 
Among the prime ideals of $\mathcal{O}_{q}(M_{m,n})$, those that are $\hc$-invariant are of particular interest. Indeed it follows from work of Goodearl and Letzter that the prime spectrum of $\mathcal{O}_{q}(M_{m,n})$ admits a stratification, into so-called $\hc$-strata, that is indexed by the set of $\hc$-invariant prime ideals of $\mathcal{O}_{q}(M_{m,n})$. We denote by $\hc$-$\spec(\mathcal{O}_{q}(M_{m,n}))$ the $\hc$-spectrum of $\mathcal{O}_{q}(M_{m,n})$; that is, the set of $\hc$-invariant prime ideals of $\mathcal{O}_{q}(M_{m,n})$. We also call these ideals $\hc$-primes. It turns out that there are only finitely many of them; as $q$ is not a root of unity, a theorem of Goodearl and Letzter \cite[Theorem II.5.12]{bgbook} show that the $\hc$-spectrum of $\mathcal{O}_{q}(M_{m,n})$ is finite. Then, using the theory of deleting-derivations, Cauchon has given a combinatorial description of $\hc$-$\spec(\mathcal{O}_{q}(M_{m,n}))$. Namely, he proved that $\hc$-$\spec(\mathcal{O}_{q}(M_{m,n}))$ is in bijection 
with the set of $m \times n$ Cauchon diagrams. We refer the reader to Cauchon's paper \cite{Cauchon} for details about the explicit correspondence between $\hc$-invariant prime ideals of $\oqmmn$ and $m\times n$ Cauchon diagrams.

An important subcollection of the $\hc$-primes in $\mathcal{O}_{q}(M_{m,n})$ is the collection of \emph{primitive} $\hc$-primes.  The second named author and Lenagan \cite{LL} and also the first two named authors and Nguyen \cite[Theorem 2.1]{BLN} showed that the primitivity of $\hc$-primes is equivalent to the corresponding Cauchon diagram having a special property which we now describe. To give this description, we require a definition.

\begin{defn} {\em
Let $C$ be an $m\times n$ labelled diagram with $k$ white squares and labels $\ell_1< \cdots < \ell_k$. We define the \emph{skew-adjacency matrix}, $M(C)$, of $C$ to be the $k\times k$ matrix whose $(i,j)$ entry is:
\begin{enumerate}
\item 1 if the square labelled $\ell_i$ is strictly to the left and in the same row as the square labelled $\ell_j$ or is strictly above and in the same column as the square labelled $\ell_j$;
\item $-1$ if the square labelled $\ell_i$ is strictly to the right and in the same row as the square labelled $\ell_j$ or is strictly below and in the same column as the square labelled $\ell_j$;
\item 0 otherwise.
\end{enumerate}
} \end{defn}

Observe that $M(C)$ is independent of the set of labels which appear in $C$. 

See, for example, Figure $3$.
\begin{figure}[h] \label{fig: 3}
\vskip 2mm
$C:~~$ \begin{tabular}{|p{0.30cm}|p{0.30cm}|p{0.30cm}|p{0.30cm}|}

\hline
$1$ & \gray & \gray &\gray 
 \\
\hline
\gray& \gray & $2$ & \gray 
 \\
\hline
$3$& \gray & $4$ & $5$  
\\
\hline
\end{tabular}\hskip 3mm $\mapsto$ \hskip 3mm $M(C) = \left( \begin{array}{rrrrr} 0&0&1&0&0\\ 0&0&0&1&0 \\ -1&0&0&1&1 \\ 0&-1&-1&0&1 \\0&0&-1&-1&0\\
\end{array} \right)$
\caption{A labelled diagram $C$ and its corresponding skew-adjacency matrix $M(C)$.}
\end{figure}

The following theorem shows that primitivity is equivalent to invertibility of the skew-adjacency matrix of the corresponding diagram.

\begin{thm} Let $P$ be an $\ch$-invariant prime in $\mathcal{O}_{q}(M_{m,n})$ and let $C$ be its corresponding Cauchon diagram. Then $P$ is primitive if and only if $\det(M(C)) \not =  0$. \label{thm: det}
\end{thm}
\noindent {\bf Proof.} See \cite[Theorem 2.1]{BLN}. \qed

With this result in mind, we define primitive diagrams.
\begin{defn} {\em Let $m$ and $n$ be natural numbers.  We say that an $m\times n$ diagram $C$ is \emph{primitive} if $\det(M(C))\not =0$.}
\end{defn}

In an earlier paper \cite{BLN}, the primitive $2\times n$ Cauchon diagrams were enumerated using Theorem \ref{thm: det}. In addition to this a conjecture was made that the determinant of a skew-adjacency matrix corresponding to a Cauchon diagram is either 0 or a power of 4. In this paper we prove a stronger version of this conjecture. 

\begin{thm} \label{thm: main3} Let $m$ and $n$ be natural numbers and let $D$ be an $m\times n$ diagram.  Then $\det( M(D))$ is either $0$ or a power of $4$.\label{thm: power4}
\end{thm}
  
We prove this result in \S \ref{Pfaffian}.
\section{Pfaffians and decompositions of diagrams}\label{Pfaffian}
In this section we give some basic facts about Pfaffians, which are the main tools we use in evaluating the determinants of skew-symmetric matrices.  We also introduce the notion of a \emph{decomposition} of a diagram, which we use to simplify the computation of the Pfaffian.

\begin{defn}
{\em Given a labelled diagram $C$, we say that
$\pi=\{\{i_1,j_1\},\ldots ,\{i_m,j_m\}\}$ is a \emph{perfect matching} of $C$ if:
\begin{enumerate}
\item{$i_1,j_1,\ldots ,i_m,j_m$ are distinct;}
\item{$\{i_1,\ldots ,i_m,j_1,\ldots ,j_m\}$ is precisely the set of labels which appear in $C$;}
\item{$i_k<j_k$ for $1\le k\le m$;}
\item{for each $k$ the white square labelled $i_k$ is either in the same row or the same column as the white square labelled $j_k$.}
\end{enumerate}
We introduce the notation
\begin{equation}
\PMC = \left\{\pi ~|~ \pi \textnormal{ is a perfect matching of } C \right\}.
\end{equation}
}
\end{defn}

For example $\left\{\left\{1,3\right\},\left\{2,8\right\},\left\{4,7\right\},\left\{5,6\right\}\right\}$ is a perfect matching on the diagram in Figure $1$.

\begin{defn}{\em
Given a perfect matching
$\pi=\{\{i_1,j_1\},\ldots ,\{i_m,j_m\}\}$ of $C$ we call the sets $\{i_k,j_k\}$ for $k=1,\ldots ,m$ the \emph{edges} of $\pi$.
We say that an edge $\{i,j\}$ of $\pi$ is \emph{vertical} if the white squares labelled $i$ and $j$ are in the same column; otherwise we say that the edge is \emph{horizontal}.}
\end{defn}
Given a perfect matching $\pi$ of $C$, we define 
\begin{equation}
{\rm sgn}(\pi) \ := \ {\rm sgn}\left(\begin{array}{ccccccc} 1 & 2 & 3 & 4 & \cdots & 2m-1 & 2m\\
i_1 & j_1 & i_2 & j_2 & \cdots & i_m & j_m \end{array} \right).\end{equation}
We note that this definition of ${\rm sgn}(\pi)$ is independent of the order of the edges (see Lovasz \cite[p. 317]{Lov}).  It is vital, however, that $i_k<j_k$ for $1\le k\le m$.  

To compute the sign of a permutation, we use \emph{inversions}.  
\begin{defn} {\em Let ${\bf x}=(i_1,i_2,\ldots ,i_n)$ be a finite sequence of real numbers.  We define ${\rm inv}({\bf x})$ to be 
$\#\{ (j,k)~|~j<k, i_j>i_k\}$.  Given another finite real sequence ${\bf y}=(j_1,\ldots ,j_m)$, we define
${\rm inv}({\bf x}|{\bf y})=\#\{(k,\ell) ~|~j_k<i_{\ell}\}$.}
\end{defn}
In fact, since ${\rm inv}({\bf x}|{\bf y})$ 
is independent of the order of elements in ${\bf x}$ and ${\bf y}$ we will also use the notation ${\rm inv}(X|Y)$ when talking about sets $X$ and $Y$ of real numbers.
 
The key fact we need is that if $\sigma$ is a permutation in $S_n$, then
\begin{equation} {\rm sgn}(\sigma) \ = \ (-1)^{{\rm inv}(\sigma(1),\ldots ,\sigma(n))}.
\end{equation}

We then define 
\begin{equation} {\rm Pfaffian}(C) \ := \ \sum_{\pi\in \PMC} {\rm sgn}(\pi).\end{equation}
In particular, if $C$ has no perfect matchings, then ${\rm Pfaffian}(C)=0$.

The significance of the Pfaffian comes from the following result (see for instance \cite[Lemma 8.2.2]{Lov}).
\begin{thm}
Let $m, n$ be natural numbers and let $C$ be an $m \times n$ labelled diagram. Then $${\rm Pfaffian}(C)^2= \det(M(C)).$$
\end{thm}
Our goal is to reduce the evaluation of the Pfaffian of a skew-symmetric matrix to counting the zeros of an associated quadratic form mod 2.

To do this we introduce the concept of a \emph{decomposition} of an $m\times n$ diagram.

\begin{defn} {\em
Let $C$ be an $m\times n$ labelled diagram with labels $\left\{1,2,\ldots,k\right\}$. We say that $(V,H)$ is a \emph{decomposition} of $C$ if:

\begin{enumerate}
\item $V\cup H = \left\{1,2,\ldots,k\right\}$ and $V\cap H =\emptyset$;
\item For $1\leq j\leq n$, $\#\left\{ i\in V~|~ \textnormal{ square }i\; \textnormal{is in column } j\right\}$ is even.
\end{enumerate}
We write \begin{equation} (V,H)\vdash C\end{equation} when $(V,H)$ is a decomposition of $C$. 
} \end{defn}
See, for example, Figure $4$.
\vskip 2mm
\begin{figure}[h] \label{fig: 5}
\vskip 2mm
 \begin{tabular}{|p{0.30cm}|p{0.30cm}|p{0.30cm}|p{0.30cm}|p{0.30cm}|}
\hline
$1$ & $2$ &$ 3$ & $4$ &\gray 
 \\
\hline
$5$ &$6$ &\gray & $7$ &$8$
 \\
\hline
\gray & $9$ & $10$ & $11$ &\gray
\\
\hline
$ 12$ &\gray & $13$ & $14$ & $15$
\\
\hline
\end{tabular}
\caption{A labelled diagram with decomposition \[V = \left\{1,5,6,8,9,10,13,15\right\},~H = \left\{2,3,4,7,11,12,14\right\}.\]}
\end{figure}

\begin{defn} {\em Let $C$ be an $m\times n$ labelled diagram and let $S$ be a subset of the labels.  We define the \emph{excess} of $S$ to be the vector \[(s_1 + 2 \mathbb{Z},s_2 + 2 \mathbb{Z},\ldots,s_m + 2 \mathbb{Z}) \in \mathbb{Z}^{m}_{2},\] 
 where  $s_i = \#\left\{ k\in S~|~ \textnormal{ square }  k \textnormal{ is in the } i\textnormal{th row of } C\right\}$.  We denote the excess of $S$ by ${\rm excess}(S)$.
When $(V,H)$ is a decomposition of $C$. We call ${\rm excess}(V)$ and ${\rm excess}(H)$ respectively the \emph{vertical-excess} and the \emph{horizontal-excess} of $(V,H)$.
} \end{defn}

The Pfaffian of a skew-symmetric matrix of a diagram can be computed in terms of decompositions of the diagram. To do this we introduce the \emph{signature} of the decomposition.  Just as we defined inversions for sequences, we introduce the notion of inversions on a subset of squares in a labelled diagram.

\begin{defn} {\em
Let $C$ be an $m\times n$ labelled diagram and let $S$ be a subset of the labels. We define the \emph{inversions} of $S$ to be \[
{\rm inv}(S)=  \sum_{i,j\in S} \delta (i<j ~\textnormal{and square }j\textnormal{ is southwest of square }i).\]
} \end{defn}

\begin{defn} {\em
Let $C$ be an $m\times n$ labelled diagram and let $(V,H)$ be a decomposition of $C$. We define the \emph{signature} of $(V,H)$ to be \[{\rm sgn}(V,H) = (-1)^{{\rm inv}(V)}\cdot(-1)^{\sum_{i<j} \delta (i\in H, j\in V)}=(-1)^{{\rm inv}(V)+{\rm inv}(V | H)}.\] 
} \end{defn}

The Pfaffian of a diagram can be expressed as a sum over decompositions.

\begin{prop}
\label{prop:pfaffian:decomposition}
Let $C$ be an $m\times n$ labelled diagram. Then \[{\rm Pfaffian}(C) =\sum{{\rm sgn}(V,H)},\] where the sum is taken over all decompositions $(V,H)$ of $C$ with ${\rm excess}(H) = (0,0,\ldots,0)$.
\end{prop}
\noindent {\bf Proof.} Recall that the Pfaffian of $C$ is the sum of the signs of the perfect matchings of $C$. We note that we can associate a decomposition $(V,H)$ with ${\rm excess}(H) = {\bf 0}$ to a perfect matching as follows. If $\left\{i_1 ,j_1 \right\},\left\{i_2 ,j_2 \right\},\ldots,\left\{i_k ,j_k \right\}$ is a perfect matching of $C$, let $V=\bigcup_\ell \left\{i_\ell , j_\ell \right\}$, where the union is over all $\ell$ such that  $i_\ell , j_\ell$ are in the same column. Similarly, $H=\bigcup_\ell \left\{i_\ell , j_\ell \right\}$, where the union runs over all $\ell$ such that $i_\ell , j_\ell$ are in the same row.
 
We note that ${\rm excess}(H)= {\bf 0}$ since $H$ is a union of pairs that are in the same row. Moreover, given any decomposition $(V,H)$ of $C$ with ${\rm excess}(H)= {\bf 0}$, we can construct a perfect matching (and possibly many).

Given a perfect matching $\pi \in \mathcal{PM}(C)$, let $(V(\pi),H(\pi))$ be the corresponding decomposition. Then 

\begin{eqnarray*}
{\rm Pfaffian}(C)&=&\sum_{\pi\in \mathcal{PM}(C)} {\rm sgn}(\pi)\\
&=& \sum_{\stackrel{  (V,H)\vdash C}{{\rm excess}(H)={\bf 0}}} \sum_{\left\{\pi\,|\, V(\pi)=V, H(\pi)=H\right\} } {\rm sgn}(\pi).
\end{eqnarray*}

To obtain the claim made in the statement of the theorem, it is sufficient to show 

\begin{equation*}
\sum_{\left\{\pi~|~ V(\pi)=V, H(\pi)=H\right\}} {\rm sgn}(\pi) = {\rm sgn}(V,H),
\end{equation*}
when ${\rm excess}(H)={\bf 0}$.
Let $V_j$ be the $m\times n$ labelled diagram in which all squares not in the $j$th column are black, and the white squares in the $j$th column are precisely the white squares in the $j$th column of $V$. We denote by $\# V_j $ the number of white squares in $V_j$. 

Similarly let $H_i$ be the $m\times n$ labelled diagram in which all squares not in the $i$th row are black and the white squares in the $i$th row are precisely the white squares in the $i$th row of $H$. We denote by $\# H_i $ the number of white squares in $H_i$.

We note that a perfect matching $\pi$ of $C$ whose corresponding decomposition is $(V,H)$ gives rise to perfect matchings $\pi_1, \pi_2 ,\ldots,\pi_n$ of $V_1 , V_2 ,\ldots,V_n$ and perfect matchings $\pi_1', \pi_2' ,\ldots,\pi_m'$ of $H_1 , H_2 ,\ldots ,H_m$, and vice versa. Moreover there are 
\begin{equation}
\left(\frac{\# V_j }{2}\right)!~\textnormal{perfect matchings of }V_j
\end{equation} and 
\begin{equation}\left(\frac{\#
H_i }{2}\right)! \textnormal{ perfect matchings of }H_i.
\end{equation} Hence there are $$\prod^{n}_{j=1} \left(\frac{\# V_j }{2}\right)! \prod^{m}_{i=1} \left(\frac{\#H_i }{2}\right)!$$ perfect matchings of $C$ whose corresponding decomposition is $(V,H)$.

Let $\pi_j$ be a perfect matching of $V_j$ for $1\leq j\leq n$ and let $\pi_i'$ be a perfect matching of $H_i$ for $1\leq i\leq m$ and let $\pi$ denote the corresponding perfect matching of $C$.

We must compute ${\rm sgn}(\pi)$. To do this, we list all pairs of edges in $\pi_1 ,\pi_2 ,\ldots, \pi_n$ followed by all pairs of edges of $\pi_1' ,\pi_2' ,\ldots, \pi_m'$ and count all the inversions. Recall that 

\begin{eqnarray*}
 {\rm sgn}(\pi)&=&(-1)^{{\rm inv}(\pi)}\\
&=& (-1)^{{\rm inv}(\pi_1 , \pi_2 ,\ldots, \pi_n ,\pi_1' ,\pi_2' ,\ldots, \pi_m')}.
\end{eqnarray*}

Note that 
\begin{eqnarray*} &~&
{\rm inv}(\pi_1 , \pi_2 ,\ldots, \pi_n ,\pi_1' ,\pi_2' ,\ldots, \pi_m')\\
& =& \sum^{n}_{j=1} {\rm inv}(\pi_j) + \sum^{m}_{i=1} {\rm inv}(\pi_i')+ \sum_{1 \leq i< j<n} {\rm inv}(\pi_i|\pi_j)\\ &~& ~~~~+ \sum^{n}_{j=1}\sum^{m}_{i=1} {\rm inv}(\pi_j|\pi_i') + \sum_{1 \leq i <  j \leq m} {\rm inv}(\pi_i'|\pi_j').
\end{eqnarray*}

Note that all labels in $\pi_i'$ are strictly less than those in $\pi_j'$ for $i<j$ since the labels of $C$ are strictly increasing form one row to the next. Hence $\sum_{1 \leq i <  j \leq m} {\rm inv}(\pi_i'|\pi_j')=0$. Note that

\begin{eqnarray*} 
\sum^{n}_{j=1} 	\sum^{m}_{i=1} {\rm inv}(\pi_j|\pi_i') 
&=&{\rm inv}(\bigcup^{n}_{j=1} \pi_j
|\bigcup^{m}_{i=1} \pi_i')\cr \\
& = & {\rm inv}(V | H).\cr \\
\end{eqnarray*}
Finally, $\sum_{1 \leq i<j \leq n} {\rm inv}(\pi_i|\pi_j) = {\rm inv}(V)$. Thus, 
\begin{equation}
{\rm inv}(\pi) = \sum^{n}_{j=1} {\rm inv}(\pi_j) + \sum^{m}_{i=1} {\rm inv}(\pi_i')+ {\rm inv}(V) + {\rm inv}(V|H).
\end{equation}
Hence \begin{equation}
{\rm sgn}(\pi) = {\rm sgn}(V,H)\prod^{n}_{j=1} {\rm sgn}(\pi_j) \prod^{m}_{i=1} {\rm sgn}(\pi_i').
\end{equation}
It follows that \begin{eqnarray*} &~&
\sum_{\stackrel{\pi\in \mathcal{PM}(C)}{V(\pi)=V, ~H(\pi)=H}} {\rm sgn}(\pi) \\
&=& {\rm sgn}(V,H)\sum_{\pi_1\in\mathcal{PM}(V_1)}\cdots \sum_{\pi_n\in\mathcal{PM}(V_n)}\sum_{\pi_1'\in\mathcal{PM}(H_1)}\cdots \sum_{\pi_m'\in\mathcal{PM}(V_m')}\prod^{n}_{j=1} {\rm sgn}(\pi_j) \prod^{m}_{i=1} {\rm sgn}(\pi_i')\cr \\
&=& {\rm sgn}(V,H)\left(\prod^{n}_{j=1} \sum_{\pi_j\in \mathcal{PM}(V_j)}{\rm sgn}(\pi_j) \right)\left(\prod^{m}_{i=1}\sum_{\pi_i'\in \mathcal{PM}(H_i)}{\rm sgn}(\pi_i')\right).
\end{eqnarray*}
We note that 
$$\sum_{\pi_j\in \mathcal{PM}(V_j)}{\rm sgn}(\pi_j)
 \ = \ \sum_{\pi_i'\in \mathcal{PM}(H_i)}{\rm sgn}(\pi_i') \ = \ 1$$ by an induction argument \cite[Lemma 2.3]{BLN}.  The result follows.  \qed

We use this result to show that ${\rm Pfaffian}(C)$ is either $0$ or $\pm 2^k$ if $C$ is a labelled diagram.  To do this we need a well-known result about the number of zeros of a quadratic form over a finite field of characteristic $2$.
\begin{lem} Let $F(y_1,\ldots ,y_d)$ be a quadratic polynomial in $\mathbb{Z}[y_1,\ldots ,y_d]$.  Then
$$\sum_{(y_1,\ldots ,y_d)\in \{0,1\}^d} (-1)^{F(y_1,\ldots ,y_k)}$$ is either $0$ or is $\pm 2^k$. \label{lem: quadform}
\end{lem}
\noindent {\bf Proof.} There exists a quadratic form $G$ such that $$F(y_1,\ldots ,y_d)=G(y_1,\ldots ,y_d) + a_0+a_1y_1+\cdots +a_d y_d.$$  Note that for $(y_1,\ldots ,y_d)\in \{0,1\}^d$,
$$F(y_1,\ldots ,y_d) \equiv G(y_1,\ldots ,y_d)+a_1y_1^2+\cdots +a_dy_d^2 + a_0 ~(\bmod \, 2).$$
Let $H$ denote the quadratic form
$$H(y_1,\ldots ,y_d):=G(y_1,\ldots ,y_d)+a_1y_1^2+\cdots +a_dy_d^2.$$
Then
$$\sum_{(y_1,\ldots ,y_d)\in \{0,1\}^d} (-1)^{F(y_1,\ldots ,y_k)}=(-1)^{a_0}\sum_{(y_1,\ldots ,y_d)\in \{0,1\}^d} (-1)^{H(y_1,\ldots ,y_k)}.$$
This sum is either $0$ or plus or minus a power of $2$ (see \cite[Theorem 16.34]{quadform}).   \qed

\begin{prop}
Let $m$ and $n$ be natural numbers and let $C$ be an $m\times n$ diagram.  Then ${\rm Pfaffian}(C)$ is either $0$ or $\pm 2^k$ for some nonnegative integer $k$.  \label{prop: pm2}
\end{prop}
\noindent {\bf Proof.} Let $N$ be the number of white squares of $C$.  We make $C$ into a labelled $m\times n$  diagram with white squares labelled $1$ through $N$.  For each white square $i$, we create a variable $x_i$ which takes on the values $0$ and $1$. For each decomposition of $C$, we let $x_i=0$ if square $i$ is in the vertical part of the matching and let $x_i=1$ if it is in the horizontal part. 
Then a decomposition $(V,H)$ of $C$ corresponds to a perfect matching if and only if it the number of squares in $H$ that are in row $i$ is even for $1\le i\le m$ and the number of squares in $V$ that are in column $j$ is even for $1\le j\le n$.   That is
\begin{equation}
\sum_{k=1}^N x_k \delta({\rm square~}k~{\rm is~in~row~}i) \equiv 0 ~(\bmod\, 2)\qquad {\rm for}~1\le i\le m
\end{equation}
and
\begin{equation}
\sum_{k=1}^N (1-x_k) \delta({\rm square~}k~{\rm is~in~column~}j) \equiv 0 ~(\bmod\, 2)\qquad {\rm for}~1\le j\le n.
\end{equation}
We define
\begin{eqnarray*}
&~& F(x_1,\ldots ,x_N) \\
& :=&  2^{-m-n} \prod_{i=1}^m\prod_{j=1}^n \left(1+(-1)^{\sum_{k=1}^N x_k \delta({\rm square~}k~{\rm is~in~row~}i)} \right) \left( 1+(-1)^{\sum_{k=1}^N (1-x_k) \delta({\rm square~}k~{\rm is~in~column~}j) } \right). 
\end{eqnarray*}
Thus a sequence $(x_1,\ldots ,x_N)\in \{0,1\}^N$ gives rise to a decomposition corresponding to a perfect matching if and only if $F(x_1,\ldots ,x_N)=1$; otherwise it is $0$.
Let $Q(x_1,\ldots ,x_N,s_1,\ldots ,s_m,t_1,\ldots ,t_m)$ denote the quadratic polynomial
\begin{equation}
\sum_{k=1}^N \sum_{i=1}^m \sum_{j=1}^n x_k s_i \delta({\rm square~}k{\rm ~is~in~row~}i) + (1-x_k)t_j\delta({\rm square~}k~{\rm is~in~column~}j).
\end{equation}
Note that 
\begin{eqnarray*}
&~& \prod_{i=1}^m\prod_{j=1}^n \left(1+(-1)^{\sum_{k=1}^N x_k \delta({\rm square~}k~{\rm is~in~row~}i)} \right) \left( 1+(-1)^{\sum_{k=1}^N (1-x_k) \delta({\rm square~}k~{\rm is~in~column~}j) } \right) \\
&=& \sum_{S\subseteq [1,m]}\sum_{T\subseteq [1,n]} (-1)^{\sum_{k=1}^N x_k \delta({\rm square~}k~{\rm is~in~a~row~indexed~by~}S)+\sum_{k=1}^N (1-x_k) \delta({\rm square~}k~{\rm is~in~column~indexed~by~}T) }\\
&=& \sum_{(s_1,\ldots,s_m)\in \{0,1\}^m} \sum_{(t_1,\ldots ,t_n)\in \{0,1\}^n} \Big[
 (-1)^{\sum_{k=1}^N\sum_{i=1}^m\sum_{j=1}^n x_k s_i \delta({\rm square~}k{\rm ~is~in~row~}i)}  \\
 &~& \hspace{5cm}  \times (-1)^{ (1-x_k)t_j\delta({\rm square~}k~{\rm is~in~column~}j)} \Big]\\
 &=& \sum_{(s_1,\ldots ,s_m,t_1,\ldots ,t_n)\in \{0,1\}^{n+ m}} (-1)^{Q(x_1,\ldots ,x_N,s_1,\ldots ,s_m,t_1,\ldots ,t_m)}.
 \end{eqnarray*}

Let $(x_1,\ldots ,x_N)\in \{0,1\}^N$ be a sequence which corresponds to a diagram corresponding to a perfect matching of $C$.  In order to compute the Pfaffian of $C$, we must compute the sign of the decomposition corresponding to this sequence.  Let 
\begin{equation}
V \ =  \ \{ i~|~x_i=0\}\qquad {\rm and}\qquad H \ = \ \{ i~|~x_i=1\}.
\end{equation}
Then 
$$(-1)^{{\rm inv}(V)} = (-1)^{\sum_{i<j} (1-x_i)(1-x_j)\delta({\rm square ~}i~{\rm is~southwest~of~square}~j )},$$
since if $x_i$ or $x_j$ is in $H$ then the terms do not contribute.
Similarly,
$$(-1)^{{\rm inv}(V|H)} = (-1)^{\sum_{i<j} x_i(1-x_j)}.$$
Since all possible decompositions of $C$ are parameterized by elements of $\{0,1\}^N$, we deduce from Proposition \ref{prop:pfaffian:decomposition} that
\begin{eqnarray*}&~& {\rm Pfaffian}(C) \\
& = & \sum_{\stackrel{  (V,H)\vdash C}{{\rm excess}(H)={\bf 0}}} {\rm sgn}(V,H)\\
& = & \sum_{\stackrel{  (V,H)\vdash C}{{\rm excess}(H)={\bf 0}}} (-1)^{{\rm inv}(V)+{\rm inv}(V,H)}\\
&=&
\sum_{(x_1,\ldots ,x_N)\in \{0,1\}^N} \Big[ (-1)^{\sum_{i<j} (1-x_i)(1-x_j)\delta({\rm square ~}i~{\rm is~southwest~of~square}~j )}(-1)^{\sum_{i<j} x_i(1-x_j)} \\
&~& \hspace{3cm}  \times ~ F(x_1,\ldots ,x_N) \Big]
\end{eqnarray*}
We note that
$$ {\sum_{i<j} (1-x_i)(1-x_j)\delta({\rm square ~}i~{\rm is~southwest~of~square}~j )}+{\sum_{i<j} x_i(1-x_j)}  $$
 is a quadratic polynomial in $x_1,\ldots ,x_N$.  We denote this polynomial by $P(x_1,\ldots ,x_N)$.
By our earlier remarks, 
$$2^{n+m} F(x_1,\ldots ,x_N) \ = \ \sum_{(s_1,\ldots ,s_m,t_1,\ldots ,t_n)\in \{0,1\}^{m+n}}
(-1)^{Q(x_1,\ldots ,x_N,s_1,\ldots ,s_m,t_1,\ldots ,t_n)},$$ for some quadratic polynomial $Q$ in 
$x_1,\ldots ,x_N,s_1,\ldots ,s_m,t_1,\ldots ,t_n$.
Thus
\begin{equation}
{\rm Pfaffian}(C) = 
2^{-n-m}\sum_{(x_1,\ldots ,x_N,s_1,\ldots ,s_m,t_1,\ldots ,t_n)\in \{0,1\}^{N+n+m} } (-1)^{P(x_1,\ldots ,x_N)+Q(x_1,\ldots ,x_N,s_1,\ldots ,s_m,t_1,\ldots ,t_n)}.
\end{equation}
Since $P+Q$ is a quadratic polynomial in $x_1,\ldots ,x_N,s_1,\ldots ,s_m,t_1,\ldots ,t_n$ and the Pfaffian of $C$ is an integer, we see that the Pfaffian of $C$ is either $0$ or $\pm 2^k$ for some natural number $k$ by Lemma \ref{lem: quadform}. \qed
\vskip 2mm
As an immediate corollary we obtain the proof of Theorem \ref{thm: power4}.
\vskip 2mm
\noindent {\bf Proof of Theorem \ref{thm: power4}.}  Let $C$ be a labelled diagram.  Then by Proposition \ref{prop: pm2} the Pfaffian of $C$ is either $0$ or $\pm 2^k$.  Since the determinant of $M(C)$ is the square of the Pfaffian of $C$ \cite[Lemma 8.2.2]{Lov}, we obtain the desired result. \qed
\vskip 2mm
A simple but useful remark can be made about the fact that nonzero Pfaffians are of the form $\pm 2^k$.
\begin{rem}  Let $C$ be a labelled diagram.  Then
$${\rm Pfaffian}(C) = 0 \textnormal{ \em if and only if }{\rm Pfaffian}(C) \equiv 0~(\bmod\, 3).$$
\label{mod3}
\end{rem}

\section{The Diagram Semigroup}\label{semigroup}

Throughout this section we fix a natural number $m$ and put a semigroup structure on the collection of all diagrams with exactly $m$ rows. This description is straightforward.

\begin{defn} {\em
Let $C_1$ and $C_2$ be respectively $m\times n_1$ and $m\times n_2$ diagrams. We define $C_1 \star C_2$ to be the $m\times (n_1+n_2)$ diagram obtained by placing $C_2$ immediately after $C_1$.
} \end{defn}

\begin{ex}
Let $m=3$ and $C_1 =\begin{tabular}{|p{0.30cm}|p{0.30cm}|p{0.30cm}|p{0.30cm}|}

\hline
 &  &  &\gray 
 \\
\hline
& \gray & & 
 \\
\hline
\gray&  & & 
\\
\hline
\end{tabular}$ and $C_2 = \begin{tabular}{|p{0.30cm}|p{0.30cm}|p{0.30cm}|}

\hline
  &  &\gray 
 \\
\hline
& \gray & 
 \\
\hline
\gray&  & 
\\
\hline
\end{tabular}$. \\ \ \\ \ \\
Then $C_1 \star C_2 = \begin{tabular}{|p{0.30cm}|p{0.30cm}|p{0.30cm}|p{0.30cm}|p{0.30cm}|p{0.30cm}|p{0.30cm}|}

\hline
 &  &  &\gray  & &  &\gray
 \\
\hline
& \gray & & & & \gray & 
 \\
\hline
\gray&  & & & \gray&  & 
\\
\hline
\end{tabular}$.
\end{ex}

The empty diagram is the identity and it is clear that $\star$ is associative.

\begin{notn}We let $\mathcal{D}_m$ denote the semigroup consisting of diagrams with exactly $m$ rows. We let $\mathcal{C}_m$ denote the subset of $\mathcal{D}_m$ consisting of Cauchon diagrams with exactly $m$ rows. \end{notn}

Our goal is to determine which diagrams in $\mathcal{D}_m$ are primitive. To do this we introduce a related object, which we call the \emph{Excess group}.

\begin{defn} {\em
Let $m$ be a natural number. We define the $m$th excess group 
\begin{equation}
{\rm Ex}_m = \left\{({\bf v},{\bf h},\varepsilon ) ~|~ {\bf v},{\bf h} \in (\mathbb{Z}/2\mathbb{Z})^m,\varepsilon \in \left\{\pm 1\right\}, {\bf v}\cdot{\bf 1} =0\right\} .
\end{equation}

with multiplication given by \begin{equation} ({\bf v},{\bf h},\varepsilon) \cdot ({\bf v'},{\bf h'},\varepsilon') = ({\bf v}+{\bf v'},{\bf h} +{\bf h'},\varepsilon'') \end{equation}
where $\varepsilon'' = \varepsilon \varepsilon' (-1)^{\sum_{1 \leq i < j \leq m} v_j (v_i' +h_i') + \sum_{1 \leq i \leq j \leq m} v_j' h_i}$, ${\bf v}=(v_1,\ldots ,v_m)$, ${\bf v'}=(v_1',\ldots ,v_m')$, ${\bf h}=(h_1,\ldots ,h_m)$, and ${\bf h'}=(h_1',\ldots ,h_m')$.
} \end{defn}

We note for further use that \begin{equation} \label{eq: order4}
({\bf v},{\bf h},\varepsilon)^4=({\bf 0},{\bf 0},+1) {\rm~~ for~ all~} ({\bf v},{\bf h},\varepsilon) \in {\rm Ex}_m.
\end{equation}

The name Excess group is motivated by associating an element of ${\rm Ex}_m$ to a decomposition $(V,H)$ to a diagram $C$ via the correspondence
\begin{equation}
(V,H) \mapsto  ({\rm excess}(V),{\rm excess}(H),{\rm sgn}(V,H)).
\end{equation}
As a notational convenience, if $(V,H)$ is a decomposition of an $m\times n$ diagram, we let
\begin{equation}
f(V,H) = ({\rm excess}(V),{\rm excess}(H),{\rm sgn}(V,H))\in {\rm Ex}_m.
\end{equation}

The multiplication rule for the excess group looks odd, but it is defined in this manner to make semigroup multiplication in $\mathcal{D}_m$ compatible with multiplication in ${\rm Ex}_m$. We make this statement more precise.
Note that if $C_1$ and $C_2$ are respectively an $m\times n_1$ and an $m\times n_2$ labelled diagram and $(V_1,H_1)$ is a decomposition of $C_1$ and $(V_2,H_2)$ is a decomposition of $C_2$, we can form a decomposition $(V_3,H_3)= (V_1,H_1) \star (V_2,H_2)$ of $C_1 \star C_2$ in a natural way, by labeling $C_1 \star C_2$ and then taking $V_3$ to be the union of the squares corresponding to $V_1$ and $V_2$ and $H_3$ to be the union of the squares corresponding to $H_1$ and $H_2$.

For example, if we take the labelled diagrams 
\vskip 2mm
$C_1=$ \begin{tabular}{|p{0.30cm}|p{0.30cm}|p{0.30cm}|p{0.30cm}|}

\hline
$1$ &  \gray& $2$ &\gray 
 \\
\hline
$3$ & $4$ & \gray & $5$
 \\
\hline
\gray& $6$ & $7$ & $8$
\\
\hline
\end{tabular}
and $C_2=$ \begin{tabular}{|p{0.30cm}|p{0.30cm}|p{0.30cm}|}

\hline
 $1$ & $2$ &\gray 
 \\
\hline
$3$ & $4$ &  \gray
 \\
\hline
\gray& \gray & $5$
\\
\hline
\end{tabular}
\vskip 2mm
with respective decompositions
\[V_1=\left\{1,3,5,8\right\}, H_1=\left\{2,4,6,7\right\}\] and \[ V_2=\left\{1,3\right\}, H_2=\left\{2,4,5\right\},\] then 
 \vskip 2mm
$C_1\star C_2=$ \begin{tabular}{|p{0.30cm}|p{0.30cm}|p{0.30cm}|p{0.30cm}|p{0.30cm}|p{0.30cm}|p{0.30cm}|}

\hline
 $1$ & \gray & $ 2$ &\gray  & $3$ & $4$ &\gray
 \\
\hline
$5$ & $6$ & \gray& $7$ & $8$ & $9$ & \gray
 \\
\hline
\gray& $10$ & $11$ & $12$ & \gray& \gray & $13$
\\
\hline
\end{tabular}
\vskip 2mm
and $(V_1,H_1)\star (V_2,H_2)$ is the decomposition \[V_3=\left\{1,3,5,7,8,12\right\}, H_3=\left\{2,4,6,9,10,11,13\right\}.\]

Note that $\left\{1,5,7,12\right\}$ are the squares in $C_1 \star C_2$ corresponding to $V_1$ and $\left\{3,8\right\}$ are the squares in $C_1 \star C_2$ corresponding to $V_2$. Similarly  $\left\{2,6,10,11\right\}$ are the squares in $C_1 \star C_2$ corresponding to $H_1$ and $\left\{4,9,13\right\}$ are the squares in $C_1 \star C_2$ corresponding to $H_2$.
\vskip 2mm

\begin{thm}
Let $(V,H)$ and $(V',H')$ respectively be decompositions of an $m\times n_1$ diagram $C$ and an $m\times n_2$ diagram $C'$. Let $(V'',H'')= (V,H) \star (V',H')$. Then 
$$f(V'',H'') \ = \   f(V,H)\cdot f(V',H')$$
holds in the semigroup ${\rm Ex}_m$.
\end{thm}
\noindent {\bf Proof.} Let $N$ and $N'$ denote the respectively the number of white squares of $C$ and $C'$. 
We may assume $C''$ has labels $\left\{1,2,\ldots,N+N'\right\}$.
We denote by $V_1$ and $V'_1$ respectively, the labels in $V''$ coming respectively from $V$ and $V'$.  Similarly, we denote by $H_1$ and $H'_1$ respectively, the labels in $H''$ coming respectively from $H$ and $H'$. So we have 
$V'' = V_1 \cup V'_1$, $H'' = H_1 \cup H'_1$.
Then \begin{equation}
{\rm sgn}(V'',H'') = (-1)^{{\rm inv}(V'')} (-1)^{\sum_{i<j} \delta(i \in H'', j\in V'')}.
\end{equation}
Note that \begin{eqnarray*} 
{\rm inv}(V'')&=&{\rm inv}( V_1 \cup V'_1)\\
&=&{\rm inv}(V_1) + {\rm inv}(V'_1) + {\rm inv}(V_1 \mid V'_1).\\
&=&{\rm inv}(V) + {\rm inv}(V') + {\rm inv}(V_1 \mid V'_1).\\
\end{eqnarray*}
Also, \begin{eqnarray*}
 \sum_{i<j} \delta(i \in H'', j\in V'')
& =& \sum_{i<j} \delta(i \in H_1, j\in V_1) + \sum_{i<j} \delta(i \in H_1, j\in V'_1)\\
 &~& ~~~~+ \sum_{i<j} \delta(i \in H'_1, j\in V_1) + \sum_{i<j} \delta(i \in H'_1, j\in V'_1).\\
 & =& \sum_{i<j} \delta(i \in H, j\in V) + \sum_{i<j} \delta(i \in H_1, j\in V'_1)\\
 &~& ~~~~+ \sum_{i<j} \delta(i \in H'_1, j\in V_1) + \sum_{i<j} \delta(i \in H', j\in V').
\end{eqnarray*}
Hence \begin{eqnarray*}
{\rm inv}(V'') + \sum_{i<j} \delta(i \in H'', j\in V'') &=& {\rm inv}(V) + \sum_{i<j} \delta(i \in H, j\in V) + {\rm inv}(V') \\
&~& ~~~+ \sum_{i<j} \delta(i \in H', j\in V') + g((V,H),(V',H')),\end{eqnarray*}
where \begin{equation}
g((V,H),(V',H')) = {\rm inv}(V_1 | V'_1) + \sum_{i<j} \delta(i \in H_1, j\in V'_1) + \sum_{i<j} \delta(i \in H'_1, j\in V_1).
\end{equation}
Consequently, \begin{equation}
{\rm sgn}(V'',H'') = {\rm sgn}(V,H) {\rm sgn}(V',H') (-1)^{g((V,H),(V',H'))}.
\end{equation}
Let \begin{eqnarray*}
{\bf v}& = &(v_1,\ldots,v_m) = {\rm excess}(V)={\rm excess}(V_1),\\
{\bf v'} &= &(v_1',\ldots,v_m') = {\rm excess}(V')={\rm excess}(V'_1),\\
{\bf h}& = &(h_1,\ldots,h_m) = {\rm excess}(H)={\rm excess}(H_1),\\
{\bf h'} &=& (h_1',\ldots,h_m') = {\rm excess}(H')={\rm excess}(H'_1).\\\end{eqnarray*}
\begin{equation}
\end{equation}
Then it is sufficient to show\begin{equation}
g((V,H),(V',H')) \equiv \sum_{1 \leq i < j \leq m} v_j (v_i' +h_i') + \sum_{1 \leq i \leq j \leq m} v_j' h_i ~(\bmod\, 2).
\end{equation}
Note that \begin{eqnarray*}
{\rm inv}(V_1 | V'_1)& =& \sum_{1\leq i<j \leq N + N'} \delta(i \in V'_1, j\in V_1) \\
&=& \sum_{1\leq i<j \leq m} \# \left\{k~|~ k\textnormal{ is in row } i {\rm ~of} V'_1\right\} \# \left\{k~|~ k \textnormal{ is in row } j \textnormal{ of } V_1\right\}\\
&\equiv & \sum_{1 \leq i<j \leq m} v_j v_i' ~(\bmod \, 2).
\end{eqnarray*}

Similarly, \begin{equation}
\sum_{i<j} \delta (i \in H_1, j \in V'_1) \equiv \sum_{1 \leq i \le j \leq m} v_j h_i'  ~(\bmod\, 2)
\label{eq: leq}
\end{equation} 
and
\begin{equation}
\sum_{i<j} \delta(i \in H'_1, j\in V_1) \equiv \sum_{1 \leq i < j \leq m} v_j h_i'  ~(\bmod\, 2).
\end{equation} 
We note that the sum on the right-hand side of equation (\ref{eq: leq}) runs over $1\le i\le j\le m$ and not $1\le i<j\le m$ because elements of $V'_1$ have greater labels than the elements of $H_1$ that are in the same row.  Putting these results together gives the desired result.    
\qed

We now show the relationship between $\mathcal{D}_m$ and ${\rm Ex}_m$.
Before doing this, we define a useful ideal of the ring $\mathbb{Z}_3 \left[{\rm Ex}_m\right]$.  We let
$J_m$ denote the ideal of $\mathbb{Z}_3 \left[{\rm Ex}_m\right]$ generated by the element $({\bf 0},{\bf 0},+1)+({\bf 0},{\bf 0},-1)$.  That is,
\begin{equation} 
J_m \ = \ \Big(({\bf 0},{\bf 0},+1)+({\bf 0},{\bf 0},-1)\Big) \ \subseteq \mathbb{Z}_3 \left[{\rm Ex}_m\right].
\end{equation}
\begin{prop}
There is a semigroup homomorphism

\begin{equation}
[\,\cdot\,] : \mathcal{D}_m \rightarrow \mathbb{Z}_3 \left[{\rm Ex}_m\right]/J_m
\end{equation}
given by
\begin{equation}
[C] = \sum_{(V,H) \vdash C} f(V,H) + J_m.
\end{equation}
Moreover if $[C_1] = [C_2]$ then $C_1$ is primitive if and only if $C_2$ is primitive.
\label{primitivemap}
\label{prop: homomorphism}
\end{prop}
\noindent {\bf Proof.} Let $C_1$ and $C_2$ be respectively  $m\times n_1$ and $m\times n_2$ diagrams. Then every decomposition of $C_1 \star C_2$ can be written as $(V_1,H_1) \star (V_2,H_2)$ for some decomposition $(V_1,H_1)$ of $C_1$ and for some decomposition $(V_2,H_2)$ of $C_2$. Moreover, the correspondence $((V_1,H_1) ,(V_2,H_2) )\mapsto (V_1,H_1) \star (V_2,H_2)$ gives a bijection between decompositions of $C_1 \star C_2$. Hence 
\begin{eqnarray*}
[C_1 \star C_2] &=& \sum_{(V,H) \vdash C_1 \star C_2} f(V,H) + J_m \\
&=&\sum_{(V_1,H_1) \vdash C_1} \sum_{(V_2,H_2) \vdash C_2}f(V_1,H_1)\cdot f(V_2,H_2)~~ + J_m \\
&=& \sum_{(V_1,H_1) \vdash C_1} f(V_1,H_1)\cdot \sum_{(V_2,H_2) \vdash C_2} f(V_2,H_2) ~~+J_m\\
&=& [C_1] \cdot [C_2].\\
\end{eqnarray*}
To show the second part of the claim, recall by Remark \ref{mod3} that ${\rm Pfaffian}(C) = 0$ if an only if ${\rm Pfaffian}(C) \equiv 0 $ (mod 3). We can write  
\begin{equation}
[C] = \sum_{({\bf v},{\bf h},\varepsilon )\in {\rm Ex}_m} a_{({\bf v},{\bf h},\varepsilon )} ({\bf v},{\bf h},\varepsilon ) + J_m,
\end{equation}
with $a_{({\bf v},{\bf h},\varepsilon )} \in \mathbb{Z}_3$.
Then the Pfaffian of $C$ mod 3 is just
\begin{equation}
\sum_{\stackrel{(V,H) \vdash C}{{\rm excess}(H)={\bf 0}}} {\rm sgn}(V,H) \equiv \sum_{ \left\{({\bf v},{\bf h},\varepsilon )\in {\rm Ex}_m: {\bf h}={\bf 0}\right\}} a_{({\bf v},{\bf h},\varepsilon )}\varepsilon~(\bmod\, 3).
\end{equation}
Thus primitivity is completely determined by the image $[C]$ of $C$ in $\mathbb{Z}_3[{\rm Ex}_m]$ and is independent of choice of representative of $[C]$ mod $J_m$.
 \qed

We note that primitivity of $C$ can be deduced purely by looking at $[C] \in \mathbb{Z}_3\left[ {\rm Ex}_m\right]/J_m$. The significance of this is that $\mathbb{Z}_3\left[ {\rm Ex}_m\right]/J_m$ is just a finite ring and so primitivity is reduced to a finite problem.  But we can say even more than this.  The proof of Proposition \ref{prop: homomorphism} shows that 
$${\rm Pfaffian}(C) \equiv \sum_{ \left\{({\bf v},{\bf h},\varepsilon )\in {\rm Ex}_m: {\bf h}={\bf 0}\right\}} a_{({\bf v},{\bf h},\varepsilon )}\varepsilon~(\bmod\, 3).$$
Let ${\bf c}={\rm Excess}(C)\in \mathbb{Z}_2^m$.  If $(V,H)$ is a decomposition of $C$ with ${\rm excess}(H)={\bf 0}$ then ${\rm excess}(V)={\bf c}$.  Hence we obtain the following remark.
\begin{rem} \label{rem:2}Let $C$ be an $m\times n$ labelled diagram and write
$$[C] \ = \ \sum_{({\bf v},{\bf h},\varepsilon )\in { \rm Ex}_m} a_{({\bf v},{\bf h},\varepsilon )} ({\bf v},{\bf h},\varepsilon ) + J_m,
$$
with $a_{({\bf v},{\bf h},\varepsilon )} \in \mathbb{Z}_3$.
  Then
$${\rm Pfaffian}(C)\equiv a_{({\bf c},{\bf 0},+1 )} - a_{({\bf c},{\bf 0},-1 )}~(\bmod\, 3),$$
where ${\bf c}\in \mathbb{Z}_2^m$ is the vector whose $ith$ coordinate is $2\mathbb{Z}$ if the number of white squares in the $ith$ row of $C$ is even and is $1+2\mathbb{Z}$ otherwise.\label{rem: 1}
\end{rem}

\section{Representation theory of the excess groups}\label{excess}
In this section we determine the degrees of the irreducible representations of ${\rm Ex}_m$ over a field whose characteristic is not $2$.  We then explicitly describe these representations in the case that $m=3$ and we are working over a field of characteristic $3$.  Fields of characteristic $3$ are especially useful to us in this context, because of Remark \ref{mod3}.
 \begin{prop} Let $m$ be a positive integer.  Then $({\rm Ex}_m,\cdot)$ is a group of order $4^m$ with identity $({\bf 0},{\bf 0},+1)$.  Furthermore ${\rm Ex}_m$ has the following properties:
 \begin{enumerate}
 \item $Z({\rm Ex}_m)=\{({\bf 0},{\bf 0},\pm 1), ({\bf 0},{\bf 1},\pm 1)\}$;
 \item if $m>1$, then the commutator subgroup ${\rm Ex}_m'$ of ${\rm Ex}_m$ is $\{ ({\bf 0},{\bf 0},\pm 1)\}$;
 \item ${\rm Ex}_m$ has $2^{2m-1}+2$ conjugacy classes;
 \item if $K$ is an algebraically closed field whose characteristic is not equal to $2$, then
 $K[{\rm Ex}_m]\cong K^{2^{2m-1}}\times M_{2^{m-1}}(K)\times M_{2^{m-1}}(K)$.
 \end{enumerate}
 \end{prop}
 \noindent {\bf Proof.} By definition of the multiplication, ${\rm Ex}_m$ is multiplicatively closed; moreover, $({\bf 0},{\bf 0},+1)$ acts as an identity element.  To check associativity is straightforward and all elements are invertible by equation (\ref{eq: order4}).  Thus ${\rm Ex}_m$ is a group.  By definition, $|{\rm Ex}_m|=4^m$.  
 To show $(1)$, note that if $({\bf v},{\bf h},+1)$ is in the center of ${\rm Ex}_m$, then
 $$({\bf v},{\bf h},+1)\cdot ({\bf 0},{\bf h'},+1)\cdot ({\bf v},{\bf h},+1)^{-1}\cdot ({\bf 0},{\bf h'},+1)^{-1} \ = \ 
 ({\bf 0},{\bf 0},(-1)^{{\bf v}\cdot {\bf h'}})=({\bf 0},{\bf 0},+1).$$
 In particular, if ${\bf v}\cdot {\bf h'}=0$ for every ${\bf h'}$, and so ${\bf v}={\bf 0}$.
 Since
 $$({\bf 0},{\bf h},+1)\cdot ({\bf v'},{\bf h'},+1)({\bf 0},{\bf h},+1)^{-1}
 ({\bf v'},{\bf h'},+1)^{-1} =({\bf 0},{\bf 0},(-1)^{{\bf v'}\cdot {\bf h}}),$$
 we see that ${\bf h}$ must be orthogonal to every vector ${\bf v}$ satisfying ${\bf v}\cdot {\bf 1}=0$.
 Since $({\bf 0},{\bf 0},-1)$ is central and the orthogonal complement of
 $$\{{\bf v}=(v_1,\ldots ,v_m) \in \mathbb{Z}_2^m ~|~v_1+\cdots +v_m=0\}$$ is spanned by ${\bf 1}$, we see that
 the center of ${\rm Ex}_m$ is 
$ \{ ({\bf 0},{\bf 0},\pm 1), ({\bf 0},{\bf 1},\pm 1)\}$, which establishes $(1)$.
 
 For $(2)$, note we have a map
 $\phi: {\rm Ex}_m\rightarrow {\mathbb Z}_2^{2m}$ which sends $ ({\bf v},{\bf h},\pm 1)$ to $({\bf v},{\bf h})$.  The image of this map is a subgroup of ${\mathbb Z}_2^{2m}$ of index $2$ and hence the kernel of this map has size $2$.  Thus ${\rm Ex}_m'\subseteq {\rm ker}(\phi)=\{ ({\bf 0},{\bf 0},\pm 1)\}$.  Since ${\rm Ex}_m$ is nonabelian for $m>1$, we see that ${\rm Ex}_m' =\{ ({\bf 0},{\bf 0},\pm 1)\}$ for $m>1$.  This establishes $(2)$.  
 
Since ${\rm Ex}_m'=\{ ({\bf 0},{\bf 0},\pm 1)\}$, we see that if $s\in {\rm Ex}_m$ then the conjugacy class of $s$ is either $\{s\}$ or $\{s\cdot ({\bf 0},{\bf 0},\pm 1)\}$; moreover, the first case occurs if and only if $s\in Z({\rm Ex}_m)$.  Thus there are
$$|Z({\rm Ex}_m)| + \frac{|{\rm Ex}_m|-|Z({\rm Ex}_m)|}{2}=4+(4^m-4)/2 = 2^{2m-1}+2$$ conjugacy classes in ${\rm Ex}_m$.  Note that if $K$ is an algebraically closed field of characteristic not equal to $2$, then
$K[{\rm Ex}_m]$ is semisimple artinian and hence isomorphic to
$$K[{\rm Ex}_m]\cong \prod_{i=1}^d M_{n_i}(K),$$ where
$d=2^{2m-1}+2$ is the number of conjugacy classes of ${\rm Ex}_m$ \cite[p. 138, Theorem 3]{reps}.  The number of $i$ such that $n_i=1$ is the order of
${\rm Ex}_m/{\rm Ex}_m'$ \cite[p. 156]{reps}, which is $2^{2m-1}$.  Thus we may assume that $n_1,n_2>1$ and $n_i=1$ for $i\ge 2$.  Since
$$\sum_{i=1}^d n_i^2=|{\rm Ex}_m|$$ (see \cite[p. 138, Corollary 2]{reps}) we have
$$n_1^2+n_2^2 + 2^{2m-1}=4^m.$$ That is, $n_1^2+n_2^2 = 2^{2m-1}$.  Moreover $n_1,n_2$ both divide $4^m$ \cite[p. 181, Proposition 5]{reps} and hence they are powers of $2$.  It follows that $n_1=n_2=2^{m-1}$. This establishes $(4)$.   \qed
\vskip 2mm 
We show as an example an explicit isomorphism between the group algebra of 
${\rm Ex}_3$ with coefficients in  $\mathbb{Z}_3$ and $\mathbb{Z}_3^{32} \times M_4 (\mathbb{Z}_3)^2$.  This isomorphism will in fact be fundamental in proving Theorem \ref{thm: main2}

$ {\rm Ex}_3 $ is generated by the six elements:
\begin{enumerate}
\item $g_1= ((1,1,0),(0,0,0),+1)$;
\item $g_2= ((0,1,1),(0,0,0),+1)$;
\item $g_3= ((0,0,0),(1,1,1),+1)$;
\item $g_4= ((0,0,0),(1,1,0),+1)$;
\item $g_5= ((0,0,0),(0,1,1),+1)$;
\item $g_6= ((0,0,0),(0,0,0),-1)$.
\end{enumerate}
We have the relations:
\begin{enumerate}
\item $g_3$ and $g_6$ are central and have order two;
\item $g_1 ^2 = g_2^2=g_6$ and $g_4^2=g_5^2=1$;
\item $[g_1, g_2]=[g_1,g_5]=[g_2,g_4]=g_6$ and $[g_1,g_4]=[g_2, g_5]=[g_4,g_5]=1$.
\end{enumerate}

Using these generators and relations, we can give an explicit isomorphism $\phi: \mathbb{Z}_3 \left[{\rm Ex}_3\right] \rightarrow (\mathbb{Z}_3)^{32} \times M_4( \mathbb{Z}_3)^2$. To give this isomorphism, it is sufficient to give the thirty-two 1-dim representations and the two 4-dim representations. The 1-dim representations are fairly straightforward.  For each vector ${\bf w}=(w_1,w_2,\ldots,w_6)\in \mathbb{Z}_2^6$, we construct a map \begin{equation} \phi_w : \mathbb{Z}_3 \left[{\rm Ex}_3\right] \rightarrow \mathbb{Z}_3\end{equation} via \begin{equation}\phi_w ({\bf v},{\bf h},\varepsilon) 
 \ = \ (-1)^{{\bf w}\cdot ({\bf v},{\bf h})}~(\bmod\, 3),\end{equation} 
where $({\bf v},{\bf h})\in \mathbb{Z}_2^6$ is the vector whose first three coordinates come from ${\bf v}$ and whose final three coordinates come from ${\bf h}$.  The operation $\cdot$ is just the ordinary dot product.  We note that although there are $64$ choices of ${\bf w}\in \mathbb{Z}_2^6$, this gives only $32$ distinct representations, since $(1,1,1)\cdot {\bf v}=0$ for every $({\bf v},{\bf h},\varepsilon)\in {\rm Ex}_3$ and hence the distinct $1$-dimensional representations correspond to cosets of $\mathbb{Z}_2^6/\langle(1,1,1,0,0,0)\rangle$.  
The two irreducible 4-dimensional representations are a little more complicated.  In describing them,  it should be understood that the entries of the matrices are taken to be in $\mathbb{Z}_3$.

In both cases the elements $g_1$ and $g_2$ can be mapped as follows:
\begin{equation}
g_1 \mapsto \left( \begin{array}{rrrr} 0&0&1&0\\ 0&0& 0&1 \\ -1&0&0&0 \\ 0&-1&0&0 \\
\end{array} \right),\hskip 2mm
g_2 \mapsto \left( \begin{array}{rrrr} 0&1&0&0\\ -1&0&0&0 \\ 0&0&0&-1 \\ 0&0&1&0 \\
\end{array} \right).\\ \end{equation}
The only generator from the list above that distinguishes these two representations is $g_3$, which is sent to the identity in one representation and is sent to minus the identity in the other.  That is,
\begin{equation}
g_3 \mapsto \pm \left( \begin{array}{rrrr} 1&0&0&0\\ 0&1&0&0 \\ 0&0&1&0 \\ 0&0&0&1 \\
\end{array} \right).\\
\end{equation}
The generators $g_4,g_5$, and $g_6$ are mapped as follows:
\begin{equation}
g_4 \mapsto \left( \begin{array}{rrrr} 1&0&0&0\\ 0&-1&0&0 \\ 0&0&1&0 \\ 0&0&0&-1 \\
\end{array} \right),\hskip 3mm
g_5 \mapsto \left( \begin{array}{rrrr} 1&0&0&0\\ 0&1&0&0 \\ 0&0&-1&0 \\ 0&0&0&-1 \\
\end{array} \right),\hskip 3mm g_6 \mapsto -\left( \begin{array}{rrrr} 1&0&0&0\\ 0&1&0&0 \\ 0&0&1&0 \\ 0&0&0&1 \\
\end{array} \right). \end{equation}
We will make use of these representations in \S \ref{enumeration}.  We note that the one-dimensional representations are killed off when we mod out by the ideal $J_3$; the reason for this is that $({\bf 0},{\bf 0},-1)$ is in the commutator and hence sent to $1$ in a $1$-dimensional representation.  But in $J_3$ the element is sent to $-1$.  The four dimensional representations, however, are unaffected by moding out by $J_3$.  We thus obtain an isomorphism
\begin{equation}
\Phi : \mathbb{Z}_3[{\rm Ex}_3]/J_3 \rightarrow M_4(\mathbb{Z}_3)\times M_4(\mathbb{Z}_3).
\end{equation}

\section{Linear recurrences of primitive $\hc$-primes}\label{automaton}

In this section we give some basic background about finite state automata and use these ideas to Prove Theorem \ref{thm: main1}. A finite state automaton is a machine that accepts as input words on a finite alphabet $\Sigma$ and has a finite number of possible outputs. We give a more formal definition. 

\begin{defn} {\em A \emph{finite state automaton} $\Gamma$ is a 5-tuple $(Q, \Sigma , \delta, q_0, F)$,where:
\begin{enumerate}
\item $Q$ is a finite set of states;
\item $\Sigma$ is a finite alphabet;
\item $\delta : Q \times \Sigma \rightarrow Q$ is a transition function;
\item $q_0 \in Q$ is the initial state;
\item $F \subseteq Q$ is the set of accepting states.
\end{enumerate} 
} \end{defn}

We refer the reader to Sipser \cite{Sipser} for more background on automata. We note that we can inductively extend the transition function $\delta$ to a function from $Q\times \Sigma^*$ to $Q$, where $\Sigma^*$ denotes the collection of finite words of $\Sigma$. We simply define $\delta(q, \varepsilon)=q$ if $\varepsilon$ is the empty word and if we have defined $\delta(q,w)$ and $x\in \Sigma$, we define $\delta(q,wx):=\delta(\delta(q,w),x)$.

\begin{defn} {\em Let $\Gamma=(Q, \Sigma , \delta, q_0, F)$ be a finite state automaton. We say that a word $w\in \Sigma^*$ is \emph{accepted} by $\Gamma$ if $\delta(q_0,w)\in F$; otherwise, we say $w$ is \emph{rejected} by $\Gamma$.}
 \end{defn}

\begin{defn} {\em
A subset $\mathcal{L}\subseteq \Sigma^*$ is called \emph{regular} if there is a finite state automaton $\Gamma$ with the property that $w\in \Sigma^*$ is accepted by $\Gamma$ if and only if $w\in \mathcal{L}$. }
\end{defn}

We take an automaton theoretic approach to the study of Cauchon Diagrams.
Let $m$ be a fixed positive integer.  Given a subset $S\subseteq \left\{1, \ldots , m\right\}$ we let $C(S)$ denote the $m\times 1$ column whose white squares appear precisely in the rows indexed by $S$. Let \begin{equation}
\Sigma_m = \left\{C(S) ~|~ S\subseteq \left\{1, \ldots , m\right\}\right\}.
\end{equation}
Then $\Sigma_m \subseteq \mathcal{D}_m$ is all $m\times 1$ diagrams and ${D}_m$ is the semigroup $\Sigma_m ^*$. We show that $\mathcal{C}_m$ is a regular subset of $\mathcal{D}_m$.

\begin{thm}
Let $m$ be a natural number. Then $\mathcal{C}_m$ is a regular subset of  $\Sigma_m ^*$.
\end{thm}
\noindent {\bf Proof.} We construct a finite state machine with input alphabet $\Sigma$ that accepts precisely the elements in  $\mathcal{C}_m$. We note that  $\mathcal{C}_m$ can be defined by the following rule:

For $2\leq i \leq m$, a column with a white square in the $i$th row cannot appear to the left of a column with a black square in the $i$th row unless the squares in this column that are in rows $1,2,\ldots,i-1$ are also colored black.

With this in mind we create accepting states $q_S$ indexed by subsets $S\subseteq \left\{2,\ldots,m\right\}$ and a single rejecting state, which we denote by $r$. Our initial state is $q_\emptyset$. The idea behind these states is that they keep track of the rows in which a white square has appeared as we move from left to right along a diagram.

Let $Q=\left\{q_S ~|~S\subseteq \left\{2,\ldots,m\right\}\right\}\cup\left\{r\right\}$. 
We define the transition $\delta: Q\times \Sigma \rightarrow Q$ as follows.  If $T\subseteq \left\{1,\ldots,m\right\}$ and 
$S\subseteq \left\{2,\ldots,m\right\}$, then $\delta (q_S , C(T))=r$ if there is an $i\in S$ such that $i \notin T$ and $j \in T$ for some $i<j$. (Note that this is saying that a white square has occurred in row $i$ somewhere and there is a column to the right of this square which has a black square in the $i$th row but a white square in the $j$th row for some $j<i$; i.e. the diagram is not Cauchon.)
Otherwise, $\delta (q_S , C(T))=q_{S\cup T \backslash \left\{1\right\}}$. Then a diagram $C\in \Sigma_m ^* = \mathcal{D}_m$ is Cauchon if and only if $\delta (q_\emptyset , C) \in \left\{q_S ~|~ S \subseteq \left\{2,\ldots,m\right\}\right\}$; i.e., if and only if $C$ is accepted by $\Gamma$. \qed

\begin{figure}[h]\begin{center}
\includegraphics{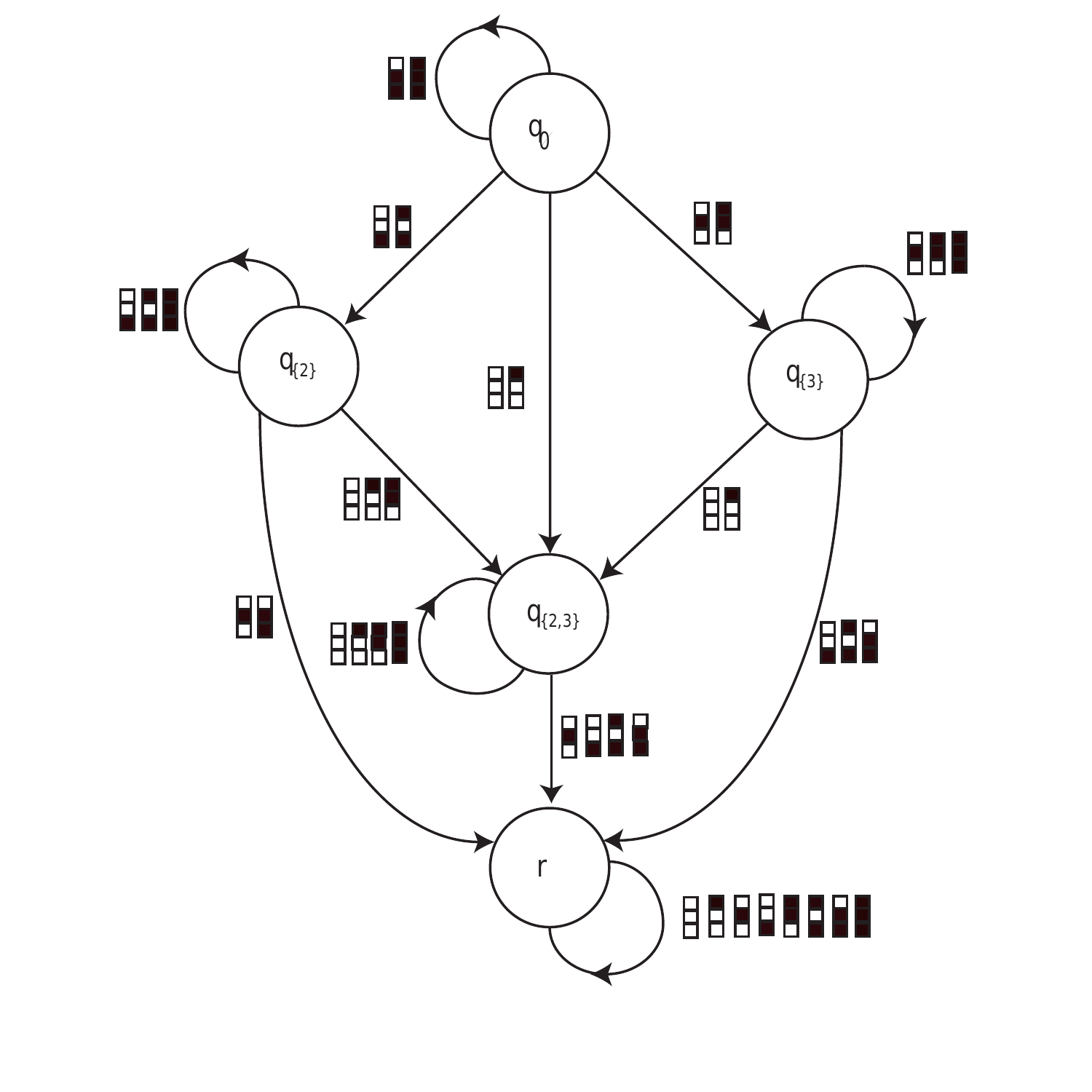}
\caption{A finite state machine that accepts $3\times n$ Cauchon diagrams.}
\end{center}
\label{fig: automaton}
\end{figure}  
The following theorem is well known, but we are unaware of a reference.

\begin{lem} Let $\Sigma$ be a finite alphabet and let $\phi : \Sigma ^* \rightarrow \mathcal{S}$ be a semigroup morphism from $\Sigma ^*$ to a finite semigroup $\mathcal{S}$. Then $\phi^{-1} (T) \subseteq \Sigma ^*$ is regular for all subsets $T\subseteq \mathcal{S}$.
\label{regular}
\end{lem}
\noindent {\bf Proof.} We let $Q=\left\{q_s ~|~ s \in \mathcal{S} \right\}$ with $q_1$ our initial state. Let $F=\left\{q_t ~|~ t \in T\right\}$. Then we have a transition function $\phi : Q\times \Sigma \rightarrow Q$ given by $\phi (q_s , x)= q_{s \phi (x)}$ for $s \in \mathcal{S}$, $x \in \Sigma$. We let $\Gamma = (Q, \Sigma , \phi, q_1, F)$. Then $w \in \Sigma ^*$ is accepted if and only if $\phi(w) \in T$. Thus $\phi^{-1}(T)\subseteq \Sigma ^*$ is a regular set. \qed

Using this result we can prove our first main result.
\vskip 2mm
\noindent {\bf Proof of Theorem \ref{thm: main1}.} We have a map $[\,\cdot\,]:\mathcal{D}_m \rightarrow \mathbb{Z}_3[{\rm Ex}_m]/J_m$. Moreover, whether or not a diagram $C\in \mathcal{D}_m$ is primitive is completely determined by its images in $\mathbb{Z}_3[{\rm Ex}_m]/J_m$. Let $T \subseteq \mathbb{Z}_3[{\rm Ex}_m]/J_m$ consist of all elements of the form $\left[C\right]$ with $C \in \mathcal{D}_m$ primitive. Then the primitive diagrams in $\mathcal{D}_m$ are given by $\left\{C\in \mathcal{D}_m~|~ \left[C\right] \in T\right\}$. By Lemma \ref{regular} the collection of primitive diagrams is a regular language in  $\mathcal{D}_m$.

Since the Cauchon diagrams are a regular sublanguage of $\mathcal{D}_m$ and an intersection of regular languages is regular \cite[Theorem 4.25]{Meduna}, we see that the primitive Cauchon diagrams form a regular language. We call this language $\mathcal{P}_m$. Then $P(m,n)$ is precisely the number of elements of the language $\mathcal{P}_m$ of length $n$. It follows that $\left\{P(m,n)\right\}^{\infty}_{n=1}$ satisfies a linear recurrence in $n$ \cite[Theorem 5.1]{Eilenberg}. \qed

\section{Enumeration of primitive $\hc$-primes in $\mathcal{O}_{q}(M_{3,n})$}\label{enumeration}

In this section we prove Theorem \ref{thm: main2}.  To do this we rely upon techniques of representation theory.  We recall that in \S \ref{excess}, we gave a description of the group algebra of the excess group ${\rm Ex}_3$ with coefficients in
  ${\mathbb Z}_3$, giving an explicit isomorphism
\begin{equation} \Phi :\mathbb{Z}_3 \left[{\rm Ex}_3\right]/J_3 \rightarrow {M}_4 (\mathbb{Z}_3)^2.
\end{equation}

At this point we look at the map $\left[\, \cdot \,\right]:\mathcal{D}_3 \rightarrow \mathbb{Z}_3 \left[{\rm Ex}_3\right]/J_3$ and classify primitive elements $\mathcal{D}_3$.
Let $G_3 \subseteq \mathbb{Z}_3 \left[ {\rm Ex}_3\right]/J_3$ be the image of $\mathcal{D}_3$ under $\left[\,\cdot\, \right]$.

\begin{prop} $G_3$ is a group of order $384$.\label{prop: 384} \end{prop}
\noindent {\bf Proof.} Since $\mathcal{D}_3$ is a semigroup and $\left[\, \cdot \, \right]:\mathcal{D}_3 \rightarrow \mathbb{Z}_3 \left[{\rm Ex}_3\right]/J_3$
is a semigroup homomorphism, $G_3$ is generated by the $8$ columns $\left\{\left[C(S)\right] ~|~ S \subseteq \left\{1,2,3\right\}\right\}$, it is sufficient to check that each of these 8 elements has an inverse.
Recall that for a subset $S\subseteq \{1,2,3\}$, we take $C(S)$ to be the $3\times 1$ column which has white squares precisely in the rows indexed by $S$; these $8$ elements are the eight columns
\vskip 2mm
\begin{tabular}{|p{0.30cm}|}
\hline
  \gray 
 \\
\hline
 \gray 
 \\
\hline
\gray
\\
\hline
\end{tabular},~~~~~~~~~~\hskip 3mm
\begin{tabular}{|p{0.30cm}|}
\hline

$1$  \\
\hline

\gray \\
\hline
\gray
\\
\hline
\end{tabular},~~~~~~~~~~\hskip 3mm
\begin{tabular}{|p{0.30cm}|}
\hline

\gray  \\
\hline

$2$ \\
\hline
\gray
\\
\hline
\end{tabular},~~~~~~~~~~\hskip 3mm\begin{tabular}{|p{0.30cm}|}
\hline

\gray  \\
\hline

\gray \\
\hline
$3$
\\
\hline
\end{tabular},~~~~~~~~~~\hskip 3mm\begin{tabular}{|p{0.30cm}|}
\hline

$1$  \\
\hline

$2$ \\
\hline
\gray
\\
\hline
\end{tabular},~~~~~~~~~~\hskip 3mm\begin{tabular}{|p{0.30cm}|}
\hline

$1$  \\
\hline

\gray \\
\hline
$3$
\\
\hline
\end{tabular},~~~~~~~~~~\hskip 3mm\begin{tabular}{|p{0.30cm}|}
\hline

 \gray \\
\hline

$2$ \\
\hline
$3$
\\
\hline
\end{tabular},~~~~~~~~~~\hskip 3mm\begin{tabular}{|p{0.30cm}|}
\hline

$1$  \\
\hline

$2$ \\
\hline
$3$
\\
\hline
\end{tabular}.
\vskip 2mm
The first four of these labeled diagrams have only one decomposition; the next three have exactly two; and the final diagram on the list has four decompositions $$(V,H)\in \left\{ (\emptyset, \{1,2,3\} ),( \{1,2 \},\{3\}),(\{1,3\},\{2\}), (\{2,3\},\{1\}) \right\}.$$
Thus the images of these eight columns under the map $[\,\cdot\,]$ are given $($mod $J_3)$ by
\begin{equation*}
[C(\emptyset)] = ({\bf 0},{\bf 0},+1),~~[C(\{1\})]=({\bf 0},(1,0,0),+1),~~[C(\{2\})]=({\bf 0},(0,1,0),+1),
\end{equation*}
\begin{equation*}
[C(\{3\})]=({\bf 0},(0,0,1),+1), ~~[C(\{1,2\})] = ((1,1,0),{\bf 0},+1)+({\bf 0},(1,1,0),+1),
\end{equation*}
\begin{equation*}
[C(\{1,3\})] = ((1,0,1),{\bf 0},+1)+({\bf 0},(1,0,1),+1),~~\end{equation*}
\begin{equation*}
[C(\{2,3\})] = ((0,1,1),{\bf 0},+1)+({\bf 0},(0,1,1),+1),~~\end{equation*}
and
\begin{eqnarray*}
[C(\{1,2,3\})] &=& ((1,1,0),(0,0,1),+1)+((1,0,1),(0,1,0),-1)\\
&~&~~~~+((0,1,1),(1,0,0),+1)+({\bf 0},(1,1,1),+1).
\end{eqnarray*}
In terms of the generators for ${\rm Ex}_3$ given in \S \ref{excess}, we see that $(\bmod\,J_3)$ we have
\begin{equation}
[C(\emptyset)] =1,~~[C(\{1\})]=g_3g_5,~~[C(\{2\})]=g_3g_4g_5,~~
\end{equation}
\begin{equation}
[C(\{3\})]=g_3g_4, ~~[C(\{1,2\})] = g_1+g_4,
\end{equation}
\begin{equation}
[C(\{1,3\})] = g_1g_2+g_4g_5,~~[C(\{2,3\})] = g_2+g_5,
\end{equation}
and
\begin{equation}
[C(\{1,2,3\})] = g_3(1+g_1g_4+g_2g_5 + g_1g_2g_4g_5).
\end{equation}

Applying the isomorphism $\Phi$ described in \S \ref{excess} to each of these elements and doing computations in matrix rings, we see that
\begin{equation}
C(\emptyset) \mapsto \left( \left( \begin{array}{rrrr} 1&0&0&0 \\ 0&1&0&0 \\ 0&0&1&0 \\ 0&0&0&1 \\
\end{array} \right), \left( \begin{array}{rrrr} 1&0&0&0 \\ 0&1&0&0 \\ 0&0&1&0 \\ 0&0&0&1 \\
\end{array}\right)\right), \end{equation}

\begin{equation}
C(\{1\}) \mapsto \left( \left(\begin{array}{rrrr} 1&0&0&0\\ 0&1&0&0 \\ 0&0&-1&0 \\ 0&0&0&-1 \\
\end{array} \right), \left(  \begin{array}{rrrr} -1&0&0&0\\ 0&-1&0&0 \\ 0&0&1&0 \\ 0&0&0&1 \\
\end{array}  \right)\right), \end{equation}

\begin{equation}
C(\{2\}) \mapsto \left( \left(  \begin{array}{rrrr} 1&0&0&0\\ 0&-1&0&0 \\ 0&0&-1&0 \\ 0&0&0&1 \\
\end{array} \right), \left(  \begin{array}{rrrr} -1&0&0&0\\ 0&1&0&0 \\ 0&0&1&0 \\ 0&0&0&-1 \\
\end{array} \right)\right), \end{equation}

\begin{equation}
C(\{3\}) \mapsto \left( \left( \begin{array}{rrrr} 1&0&0&0 \\ 0&-1&0&0 \\ 0&0&1&0 \\ 0&0&0&-1 \\
\end{array} \right), \left( \begin{array}{rrrr} -1&0&0&0 \\ 0&1&0&0 \\ 0&0&-1&0 \\ 0&0&0&1 \\
\end{array}\right)\right), \end{equation}

\begin{equation}
C(\{1,2\}) \mapsto \left( \left( \begin{array}{rrrr} 1&0&1&0 \\ 0&-1&0&1 \\ -1&0&1&0 \\ 0&-1&0&-1 \\
\end{array} \right), \left( \begin{array}{rrrr} 1&0&1&0 \\ 0&-1&0&1 \\ -1&0&1&0 \\ 0&-1&0&-1 \\
\end{array}\right)\right), \end{equation}

\begin{equation}
C(\{1,3\}) \mapsto \left( \left( \begin{array}{rrrr} 1&0&0&-1 \\ 0&-1&1&0 \\ 0&-1&-1&0 \\ 1&0&0&1 \\
\end{array} \right), \left( \begin{array}{rrrr} 1&0&0&-1 \\ 0&-1&1&0 \\ 0&-1&-1&0 \\ 1&0&0&1 \\
\end{array}\right)\right), \end{equation}

\begin{equation}
C(\{2,3\}) \mapsto \left( \left( \begin{array}{rrrr} 1&1&0&0 \\ -1&1&0&0 \\ 0&0&-1&-1 \\ 0&0&1&-1 
\end{array} \right), \left( \begin{array}{rrrr} 1&1&0&0 \\ -1&1&0&0 \\ 0&0&-1&-1 \\ 0&0&1&-1 
\end{array}\right)\right), \end{equation}

\begin{equation}
C(\{1,2,3\}) \mapsto \left( \left( \begin{array}{rrrr} 1&1&1&1 \\ -1&1&1&-1 \\ -1&-1&1&1 \\ -1&1&-1&1 \\
\end{array} \right), \left( \begin{array}{rrrr} 1&1&1&1 \\ -1&1&1&-1 \\ -1&-1&1&1 \\ -1&1&-1&1 \\
\end{array}\right)\right), \end{equation}
where the entries of the matrices are in $\mathbb{Z}_3$.
Since each of these matrices is invertible, $G_3$ is isomorphic to a subgroup of ${\rm GL}_4(\mathbb{Z}_3)^2$.  Using the computer algebra package SAGE, which can compute orders of linear groups mod $d$, we find that the group $G_3$ has order $384$. \qed
\begin{rem} We note that the image of elements of $G_3$ in ${\rm GL}_4(\mathbb{Z}_3)^2$ described in the proof of Proposition \ref{prop: 384} are always of the form $(A,A)$ or $(A,-A)$ for some invertible matrix $A$; moreover, the elements of the form $(A,A)$ are precisely those elements in $G_3$ which correspond to elements of $\mathcal{D}_3$ with an even number of white squares. \label{rem: 384}
\end{rem}
 \begin{prop} Let $C_{1} $ and $C_{2} $ be two  diagrams. Then the  determinant of the skew-adjacency matrix corresponding to $C_{1} \star C_{2} $ is the same as the determinant of the skew-adjacency matrix corresponding to $C_{2} \star C_{1} $.

\label{THM:normal} \end{prop}
\noindent {\bf Proof. } We note that we can assign labels $\{1,2,\ldots ,N\}$ to the white squares of a diagram $C$ with $N$ white squares (in such a way that $C$ is not necessarily what we call a labelled diagram) and then create an $N\times N$ skew-symmetric matrix whose $(i,j)$ entry is:
\begin{enumerate}
\item  $+1$ if the square labelled $i$ is either in the same row and strictly to the left of the square labelled $j$, or is in the same column and strictly above the square labelled $j$; 
\item $-1$ if the square labelled $i$ is either in the same row and strictly to the right of the square labelled $j$, or is in the same column and strictly below the square labelled $j$; 
\item $0$ otherwise.
\end{enumerate}
If we do this, the resulting skew-symmetric matrix is similar to the skew-adjacency matrix of $C$ in which the matrix giving the similarity is a permutation matrix.

Suppose $C_1$ and $C_2$ have $N_1$ and $N_2$ white squares respectively.  We assign the labels $\{1,2,\ldots ,N_1\}$ to $C_1$ by declaring that a white square labelled $i$ is in a column that is strictly to the left of a white square labelled $j$, then $i<j$ and if a white square labelled $i$ is in the same column and strictly above a square labelled $j$ then $i<j$.  In an analogous manner, we assign the  labels $\{1,2,\ldots ,N_2\}$ to the white squares of $C_2$; and the labels $\{1,2,\ldots ,N_1+N_2\}$ to the white squares of $C_1\star C_2$ and $C_2\star C_1$.  
Let $M_{1} $ and $M_{2} $ denote respectively the skew-symmetric matrices described above corresponding to these labellings of $C_{1} $ and $C_{2} $. The skew-symmetric matrix corresponding to $C_{1} \star C_{2} $ with this labelling is equal to \[\left( \begin{array}{cc} M_{1} &  B \\
 -B^{ T } & M_{2}  \end{array}\right).\]  

Similarly, the skew-symmetric matrix  corresponding to $C_{2} \star C_{1}$ is given by

\[\left( \begin{array}{cc} M_{2} & B^{T}\\ -B & M_{1} \end{array}\right).\]

Observe that if we take this matrix, considering it as a block $2\times 2 $ matrix,  and interchange the two rows and then interchange the two columns, we obtain the new matrix

\[\left(\begin{array}{cc} M_{1} & -B\\ B^{T} & M_{2} \end{array}  \right).\]

We again think of this matrix as being ${2 \times }{2}$ block matrix, and now multiply the first column by $-1$ and multiply the first row by $-1 $. Doing this, we obtain the matrix corresponding to the diagram $C_{1} \star C_{2}$. The result now follows. Since the skew-adjacency matrices of $C_1\star C_2$ and $C_2\star C_1$ are both similar to this matrix, we see they are either both invertible or they both fail to be invertible.  The result follows. \qed
\vskip 2mm
As an immediate corollary, we obtain the following fact.

\begin{cor} The collection of elements in the group $G_{3} $ corresponding to primitive $3\times n$ diagrams is a union of conjugacy classes.
\label{conjugacyresult}
  \end{cor}
\noindent {\bf Proof. } Let $C$ be a primitive $3\times n $ diagram, and suppose that $D$ is another $3\times n $ diagram. Then there exists an $3\times n $ diagram $D' $ such that the image of  $D' $ in $G_{ 3 } $ is the inverse of the image of $D$. Then $[D ] \cdot [C] \cdot [D']$ is conjugate to $[C ] $ in $G_3$. Hence it is sufficient to show that the diagram $D\star C \star D'$ is primitive. By Proposition \ref{THM:normal}, the determinant of the skew-adjacency matrix corresponding to $D\star C \star D' $ is the same up to sign as the determinant of the skew-adjacency matrix corresponding to $D' \star D\star C$. Since the image of $D' \star D $ in $G_3$ is the identity and the determinant of the skew-adjacency matrix corresponding to $C$ is nonzero, we see that the diagram $D\star C \star  D'$ is also primitive. The result follows. \qed
\vskip 2mm
We make a further simplification to reduce the size of $G_3$.   Let 
\begin{equation} N = \left\{[D]~|~ D \in \mathcal{D}_3
 ~{\rm and}~ [ D]= \pm ({\bf v},{\bf 0},\varepsilon) ~{\rm for~ some } ~
  {\bf v} \in \mathbb{Z}^{3}_{2},~ \varepsilon
   \in \left\{\pm 1 \right\}\right\}. \end{equation}
Clearly, $N$ is a subgroup of $G_3$ and $|N| \leq 8$. A simple computation shows that $g_1g_6+J_3=[C(\{1\})][C(\{2\})][C(\{1,2\})]^2$, 
$g_2g_6+J_3=[C(\{2\})][C(\{3\})][C(\{2,3\})]^2$, and $g_6+J_3=[C(\{1,2\})]^4$. Hence $g_1+J_3$, $g_2+J_3$ and $g_6+J_3$ belong to $N$, and we easily deduce from this that 
$N =\{1, g_6+J_3, g_1+J_3, g_1g_6 +J_3,  g_2+J_3, g_2g_6+J_3 , g_1g_2+J_3, g_1g_2g_6+J_3 \}$. In particular, $N$ is generated by  $g_1+J_3$, $g_2+J_3$ and $g_6+J_3$. Using this, it is easy to prove that $N$ is a normal subgroup of $G_3$.

\begin{prop} If $D_1,D_2 \in \mathcal{D}_3$ and $\left[D_1\right]N= \left[D_2\right]N$ then either $D_1,D_2$ are both primitive or both fail to be primitive. 
\label{prop: norm}
\end{prop}
\noindent {\bf Proof.} We put an equivalence relation $\sim$ on $\mathcal{D}_3$ by declaring $D_1 \sim D_2$ if for each $D \in \mathcal{D}_3$, $D_1 \star D$ is primitive if and only if $D_2 \star D$ is primitive. Let $$\widetilde{N}= \left\{\left[D\right]\in G_3 ~|~ D \sim C(\emptyset)\right\},$$
i.e., $\widetilde{N}$ consists of all elements equivalent to the identity.

 We claim $\widetilde{N}$ is a normal subgroup of $G_3$. To see this, note that if $\left[D_1\right],\left[D_2\right]\in \widetilde{N}$ then $D_1 \star D_2 \star D$ is primitive if and only if $D_2 \star D$ is primitive since $D_1 \sim C(\emptyset)$. But $D_2 \star D $ is primitive if and only if $D$ is primitive since $D_2 \sim C(\emptyset)$. Hence $\left[D_1 \star D_2\right]= \left[D_1\right]\left[D_2\right] \in \widetilde{N}$. Thus $\widetilde{N}$ is a group. 
 
 To see $\widetilde{N}$ is normal,  let $\left[D_0\right]\in \widetilde{N}$ and let $\left[E\right]\in G_3$. Then there exists $F\in \mathcal{D}_3$ such that $\left[E\right]\left[F\right]=1$. It is sufficient to show $\left[E\right]\left[D_0\right]\left[F\right]\in \widetilde{N}$. Note that for $D\in \mathcal{D}_3$, $E\star D_0 \star F \star D$ is primitive if and only if $D_0\star F \star D \star E$ is primitive by Proposition \ref{THM:normal}. Since $D_0\in \widetilde{N}$, $D_0\star F \star D \star E$ is primitive if and only if $F \star D \star E$ is primitive; but $F \star D \star E$ is primitive if and only if $E \star F \star D$ is primitive by Propostion \ref{THM:normal}. Since $\left[E \star F \star D\right]=\left[E\right]\left[F\right]\left[D\right]=\left[D\right]$, we see this occurs if and only if $D$ is primitive; hence $\widetilde{N}$ is normal.
 
We now show $N\subseteq \widetilde{N}$. Note that if $\left[D\right]\in N$, then $\left[D\right]=  ({\bf d},{\bf 0},\varepsilon)+J_3$. Then if $\left[C\right]= \sum a_{({\bf v},{\bf h},\varepsilon)} ({\bf h},{\bf 0},\varepsilon)+J_3$, it follows from Remark \ref{rem:2} that $C$ is primitive if and only if $$a_{({\bf c},{\bf 0},+1)}-a_{({\bf c},{\bf 0},-1)}\not\equiv 0 ~(\bmod\, 3),$$ where ${\bf c}={\rm excess}(C)$.  Since $({\bf 0},{\bf 0},+1)\equiv -({\bf 0}, {\bf 0},-1) ~(\bmod\, J_3)$, we can choose a representative of $[C]$ mod $J_3$ of the form
$$[C] \ = \sum b_{{\bf v},{\bf h}} ({\bf v},{\bf h},+1)+J_3$$ with $b_{{\bf v},{\bf h}}\in \mathbb{Z}_3$.  Primitivity of $C$ is then equivalent to $b_{{\bf c},{\bf 0}}\not \equiv 0 ~(\bmod\, 3)$.  But $\left[D\right]\left[C\right]$ is then
$$\sum   b_{{\bf v},{\bf h}} ({\bf v}+{\bf d},{\bf h},\varepsilon_{{\bf v},{\bf h}})+J_3,$$ with $\varepsilon_{{\bf v},{\bf h}} \in \{\pm 1\}$.  Thus we deduce from Remark \ref{rem:2} that $D\star C$ is primitive if and only if $\varepsilon_{{\bf c},{\bf 0}} b_{{\bf c},{\bf 0}}\not\equiv 0 ~(\bmod \, 3)$.  Since $\varepsilon_{{\bf c},{\bf 0}}$ is a unit mod $3$, we see that $D\star C$ is primitive if and only if $C$ is primitive and so $[D]\in \widetilde{N}$.  The result follows. \qed
\vskip 2mm
Thus we can look at the image of diagrams in $G_3/N$ instead of in $G_3$ to determine whether they are primitive or not.  Since $N$ has order $8$ and $G_3$ has order 384, we see that $G_3/N$ is a group of order 48; in fact, we can describe this group very well.
\begin{prop} We have an isomorphism $$\psi :G_3/N\rightarrow S_4\times \{\pm 1\}$$ given by
\[\begin{array}{lr}
\left[ C( \left\{1,2,3\right\})\right] ~ \mapsto ~  ((124),-1), &
 \left[ C( \left\{1,2\right\})\right] ~ \mapsto ~ ((1243),1),\\
 \left[ C( \left\{1,3\right\})\right]~ \mapsto ~ ((1324),1),&
 \left[ C( \left\{2,3\right\})\right]~ \mapsto ~((1234),1),\\
 \left[ C( \left\{1\right\})\right] ~ \mapsto ~ ((13)(24),-1),&
\left[ C( \left\{2\right\})\right] ~ \mapsto ~ ((12)(34),-1),\\
 \left[ C( \left\{3\right\})\right] ~ \mapsto ~ ((14)(23),-1),&
 \left[ C( \emptyset)\right]~ \mapsto ~ ({\rm id},1).\end{array}\]
 Moreover, the elements of $G_3/N$ in the preimage $\psi^{-1}((S_4\times \{1\})$ are precisely the images of the diagrams $C\in \mathcal{D}_3$ with an even number of white squares.
\end{prop}
\noindent {\bf Proof.}
Recall that we have a map $\left[\, \cdot\, \right]: \mathcal{D}_3 \rightarrow \mathbb{Z}_3 \left[ {\rm Ex}_3\right]/J_3$ in which the image of $ \mathcal{D}_3$ is a group $G_3$.  Thus $\left[\, \cdot\, \right]$ induces a surjective semigroup homomorphism from $ \mathcal{D}_3$ to $G_3$.  Doing straightforward computations with the generators and relations for $G_3$ given in \S 5, we find that the map $\psi$ defined in the statement of the proposition is indeed an isomorphism.  One sees that the generators of $\mathcal{D}_3$ that are sent to elements of the form $(\sigma, 1)$ are precisely those generators with an even number of white squares (see also Remark \ref{rem: 384}).  The result follows. \qed
\vskip 2mm
Let \begin{equation} h:G_3\rightarrow G_3/N \end{equation}
be the canonical surjection and let \begin{equation}
H:\mathcal{D}_3\rightarrow S_4\times \{\pm 1\}\end{equation}
 be given by the composition \begin{equation}
H \ = \  \psi \circ h\circ [\, \cdot \, ].\end{equation}
\begin{prop}
The primitive diagrams in $\mathcal{D}_3$ are the preimage under $H$ of the three conjugacy classes with representatives $({\rm id},1)$, $((1,2,3),1)$, and $((1,2,3,4),1)$. 
\end{prop}
\noindent {\bf Proof.} By Theorem \ref{primitivemap} and Proposition \ref{prop: norm}, primitivity of a diagram $C\in \mathcal{D}_3$ can be deduced by looking at its image under $H$.  Thus there is a subset $S \subseteq S_4 \times \left\{ \pm 1\right\}$ such that the set of primitive elements of $\mathcal{D}_3$ is precisely the preimage of $S$ under $H$. That is, \begin{equation}
\left\{C \in \mathcal{D}_3 ~|~ C\textnormal{ is primitive }\right\}= \bigcup_{s\in S} H^{-1}(\left\{s\right\}). \end{equation}
By Corollary~\ref{conjugacyresult} and Proposition \ref{prop: norm}, $S$ is a union of conjugacy classes. Note that $S_4 \times \left\{ \pm 1\right\}$ has $10$ conjugacy classes, and it is sufficient to pick a representative from each one and check the primitivity of an element in $\mathcal{D}_3$ whose image under $H$ is this representative. We note that any conjugacy class containing $(\sigma,-1)$ for some $\sigma \in S_4$ has the property that $\psi ( (\sigma,-1))$ consists entirely of diagrams with an odd number of white squares. Since an $m\times m$ skew-symmetric matrix with $m$ odd is not invertible, we see that none of these conjugacy classes correspond to primitive diagrams. This leaves 5 conjugacy classes to check. Using the isomorphism given above, we find 

\begin{equation}
\psi ( \left[ C( \emptyset)\right])  \\ =  \\ ({\rm id},1)\end{equation}
\begin{equation}
\psi ( \left[ C( \left\{1,2\right\})\right])  \\ = \\  ((1243),1)\end{equation}
\begin{equation}
\psi ( \left[ C( \left\{1\right\})\star C(\left\{2\right\})\right]) \\ =  \\ ((14)(23),1)\end{equation}
\begin{equation}
\psi ( \left[ C( \left\{1,2,3\right\})^{\star \,2} \right])  \\ = \\  ((124),1)\end{equation}
\begin{equation}
\psi ( \left[ C( \left\{1,2,3\right\})\star C( \left\{2,3\right\})\star C( \left\{1,2,3\right\}) \right])  \\ = \\  ((13),1).\end{equation}
These five elements correspond respectively to the five Cauchon diagrams
\vskip 2mm
\begin{tabular}{|p{0.30cm}|}
\hline
  \gray 
 \\
\hline
 \gray 
 \\
\hline
\gray
\\
\hline
\end{tabular},~~~~~~~~~~\hskip 3mm
\begin{tabular}{|p{0.30cm}|}
\hline

 \\
\hline

 \\
\hline
\gray
\\
\hline
\end{tabular},~~~~~~~~~~\hskip 3mm
\begin{tabular}{|p{0.30cm}|p{0.30cm}|}
\hline
  & \gray 
 \\
\hline
 \gray &
 \\
\hline
\gray & \gray
\\
\hline
\end{tabular},~~~~~~~~~~\hskip 3mm
 \begin{tabular}{|p{0.30cm}|p{0.30cm}|}
\hline
  &
 \\
\hline
 &
 \\
\hline
 & 
\\
\hline
\end{tabular},~~~~~~~~~~\hskip 3mm
\begin{tabular}{|p{0.30cm}|p{0.30cm}|p{0.30cm}|}
\hline
  &\gray  &
 \\
\hline
&  &
 \\
\hline
 & &
\\
\hline
\end{tabular}.
\vskip 2mm
We compute the corresponding skew-symmetric matrices for each of these five diagrams and find the first, second, and fourth are primitive diagrams and the third and fifth are not.  \qed
\vskip 2mm  
We now explain how we will use these techniques to enumerate the primitive $3\times n$ Cauchon diagrams. Note that the semigroup morphism $H:\mathcal{D}_3\rightarrow S_4\times \{\pm 1\}$ extends to a map $\widehat{H} : \mathbb{C}\left[\mathcal{D}_3\right]\left[\left[t\right]\right]\rightarrow \mathbb{C}\left[S_4 \times \left\{\pm1\right\}\right]\left[\left[t\right]\right]$.
The following proposition gives an expression for the generating series of the Cauchon diagrams in $\mathbb{C}\left[\mathcal{D}_3\right][[t]]$ as a rational function. To give this expression, we introduce the following functions. We let \begin{equation}
F_1 (t) = \frac{1}{1-( C(\emptyset)+C(\left\{1\right\}) )t},
\end{equation}
\begin{equation}
F_2 (t) = F_1 (t)  \star \Big[C(\left\{2\right\})+C(\left\{1,2\right\})\Big]t \star \frac{1}{1-(C(\emptyset)+C(\left\{2\right\})+C(\left\{1,2\right\}))t},
\end{equation}
\begin{equation}
F_3 (t) = F_1 (t) \star \Big[C(\left\{3\right\})+C(\left\{1,3\right\})\Big]t \star \frac{1}{1-(C(\emptyset)+C(\left\{3\right\})+C(\left\{1,3\right\}))t},
\end{equation}
and \begin{equation}
F_4 (t) = \frac{1}{1-(C(\emptyset)+C(\left\{3\right\})+C(\left\{2,3\right\})+C(\left\{1,2,3\right\}))t}.
\end{equation}
\begin{prop}\label{prop: C}
Let $${\sf C}(t) = \sum_{C\in \mathcal{C}_3} C\cdot t^{ \# {\, \rm columns \,of\,} C} \in \mathbb{C}\left[\mathcal{D}_3\right]\left[\left[t\right]\right]$$ be the generating function for Cauchon diagrams. Then 
\begin{eqnarray*} {\sf C}(t) &=& F_1(t) \star \Big[1+ (C(\left\{2,3\right\})+C(\left\{1,2,3\right\}))t \star F_4 (t)\Big] \\
&~&~~~~~+~~ F_2(t) \star \Big[1+ (C(\left\{3\right\})+ C(\left\{2,3\right\})+C(\left\{1,2,3\right\}))t \star F_4 (t)\Big]\\
&~&~~~~+ ~~
 F_3(t) \star \Big[1+ ( C(\left\{2,3\right\})+C(\left\{1,2,3\right\}))t \star F_4 (t)\Big].
 \end{eqnarray*}
\end{prop}
\noindent {\bf Proof.} Note that $\mathcal{C}_3$ consists of all elements in $\Sigma^*$ that are accepted by the automaton in Figure~\ref{fig: automaton}. To be accepted, a word must be sent to one of the states $q_{\emptyset},q_{\left\{2\right\}},q_{\left\{3\right\}},q_{\left\{2,3\right\}}$. The generating function for words that are sent to $q_{\emptyset}$ is given by $$ \sum^{\infty}_{n=0} [ C(\emptyset) + C(\left\{1\right\})]^n t^n = \frac{1}{1-[C(\emptyset)+C(\left\{1\right\})]t},$$ which is $F_1(t)$.  
The generating function for words that are sent to $q_{\left\{2\right\}}$ is given by $$  \frac{1}{1-[C(\emptyset)+C(\left\{1\right\})]t} \star \Big[C(\left\{2\right\})+C(\left\{1,2\right\})\Big]t \star \frac{1}{1-[C(\emptyset)+C(\left\{2\right\})+C(\left\{1,2\right\})]t},$$ which is $F_2(t)$.  Similarly the generating function for words that are sent to $q_{\left\{3\right\}}$ is given by $$  \frac{1}{1-[C(\emptyset)+C(\left\{1\right\})]t} \star \Big[C(\left\{3\right\})+C(\left\{1,3\right\})\Big]t \star \frac{1}{1-[C(\emptyset)+C(\left\{3\right\})+C(\left\{1,3\right\})]t},$$
which is $F_3(t)$.  The most complicated component of the generating function to count is the words in $\Sigma^*$ that are sent to $q_{\left\{2,3\right\}}$ since there are multiple paths in the automaton. This is given by 
\begin{eqnarray*}
 &~& F_1(t) \star \Big[(C(\left\{2,3\right\})+C(\left\{1,2,3\right\}))t \star F_4 (t)\Big] \\
&~&~~~~~+~~ F_2(t) \star \Big[(C(\left\{3\right\})+ C(\left\{2,3\right\})+C(\left\{1,2,3\right\}))t \star F_4 (t)\Big]\\
&~&~~~~+ ~~
 F_3(t) \star \Big[ ( C(\left\{2,3\right\})+C(\left\{1,2,3\right\}))t \star F_4 (t)\Big].
 \end{eqnarray*}
Putting these results together, we obtain the desired result. \qed

This result, while complicated, gives a way of expressing the generating function for $\left\{P(3,n)\right\}$, the number of primitive $\hc$-primes in $\mathcal{O}_{q}(M_{m,n})$, as a rational power series in $t$. 
We let \begin{equation}
G_1 (t) = \widehat{ H } (F_1 (t)) =  \frac{1}{1-[ ({\rm id},1) +((13)(24),-1)]t},
\end{equation}
\begin{equation}
G_2 (t) = \widehat{ H } (F_2 (t)), \end{equation}
that is,
\begin{equation*}
G_2(t)=G_1 (t)  \cdot \Big[  ((12)(34),-1)+((1243),1)\Big]t \cdot \frac{1}{1-[({\rm id},1)+((12)(34),-1)+((1243),1)]t}.\end{equation*}
We let
\begin{equation}
G_3 (t) = \widehat{ H } (F_3 (t)), \end{equation}
that is
\begin{equation*}
G_3 (t) = G_1 (t) \cdot \Big[ ( (14)(23),-1)+(( 1324),1)  \Big]t \cdot \frac{1}{1-[    ({\rm id},1) +((14)(23),-1)+((1324),1)  ]t},
\end{equation*}
and \begin{equation}
G_4 (t) = \widehat{ H } (F_4 (t)) =    \frac{1}{1- [ ({\rm id},1)+((14)(23),-1)+((1234),1)+((124),-1) ]t}. \end{equation}
Using these facts along with Proposition \ref{prop: C}, we get the following result.
\begin{rem} \label{lem: psiC}
We have
\begin{eqnarray*}\widehat{ H } ({\sf C}(t))  &=& G_1(t) \cdot \Big[1+ \Big(  ((1234),1)+ ((124) ,-1)  \Big)t \cdot G_4 (t)\Big] \\
&~&~~~~~+~~ G_2(t) \cdot \Big[1+ \Big(    ((14)(23),-1)+ ((1234),1)+((124),-1) \Big)t \cdot G_4 (t)\Big]\\
&~&~~~~+ ~~
 G_3(t) \cdot \Big[1+ \Big( ( (1234),1) +((124),-1)   \Big)t \cdot G_4 (t)\Big].
 \end{eqnarray*}
 \end{rem}
Using this result, we can enumerate the $3\times n$ primitive Cauchon diagrams.
\vskip 2mm
\noindent {\bf Proof of Theorem \ref{thm: main2}.}
Let $$P: S_4 \times \left\{\pm 1\right\} \rightarrow \mathbb{C}$$ be the class function that sends the conjugacy classes containing $$({\rm id},+1),~~~(
(123),+1),~~~((1234),+1)$$ to 1 and all other conjugacy classes to $0$. Then $P$ is a linear combination of characters. Using the orthogonality relations, we find
\begin{eqnarray*}
P &=& \frac{15}{48}  (\chi_1\otimes \Delta_0 + \chi_1 \otimes \Delta_1 )+ \frac{3}{48}  (\chi_2\otimes \Delta_0 +\chi_2 \otimes \Delta_1 )\\ &~&~~~~ -\frac{6}{48} (\chi_3 \otimes \Delta_0 + \chi_ 3 \otimes \Delta_1 )-\frac{3}{48} (\chi_4 \otimes \Delta_0+ \chi_ 4 \otimes \Delta_1 ) +\frac{9}{48} (\chi_ 5\otimes \Delta_0 + \chi_ 5 \otimes \Delta_1 ),
\end{eqnarray*}
where the characters $\chi_1, \chi_2, \chi_3, \chi_4, \chi_5$ are the characters of $S_4$ from the character table given in Figure $6$; $\Delta_0, \Delta_1$ are the characters of $\{\pm 1\}$ given in Figure $6$; and $\chi_i \otimes \Delta_j((\sigma,\varepsilon)) = \chi(\sigma)\Delta(\varepsilon)$.
\vskip 2mm
\begin{figure}[h]
\label{fig: 12}
\vskip 2mm
\begin{tabular}{|c|r|r|r|r|r|}
\hline
~ & ~ &~ & ~ &~ &  \\
~~ &$ \{{\rm id}\} $ & $\{(12)\}$ &$ \{(12)(34)\}$ & $\{(123)\}$ &$\{(1234)\}$ \\
~ & ~ &~ & ~ &~ &  \\
\hline
$\chi_1$ & $1$ & $1$ & $1$ & $1$ & $1$ \\
\hline
 $\chi_2$ & $1$ & 
 $-1$  & $1$ & $1$ & $-1$ \\
\hline
$ \chi_3$ & $2$ & $0$  & $2$ & $-1$ & $0$\\ \hline
$ \chi_4$ & $3$ & $1$  & $-1$ & $0$ & $-1$\\ \hline
$ \chi_5$ & $3$ & $-1$  & $-1$ & $0$ & $1$\\ \hline
\end{tabular}
\vskip 5mm 
\begin{tabular}{|c|r|r|}
\hline
~ & ~ &~  \\
~~ &$ \{ +1 \} $ & $\{-1 \}$ \\
~ & ~ & ~ \\
\hline
$\Delta_0$ & $1$ & $1$ \\
\hline
 $\Delta_1$ & $1$ & 
 $-1$   \\
\hline
\end{tabular}
\label{fig: diagramxxx}
\caption{The character tables of $S_4$ and $\{\pm 1\}$}
\end{figure}

Then the class functions $$P, \chi_1\otimes \Delta_0,\ldots ,\chi_5\otimes \Delta_0, \chi_1\otimes \Delta_1,\ldots ,\chi_5\otimes \Delta_1$$ can each be extended to functions on $\mathbb{C}\left[S_4 \times \left\{ \pm 1\right\}\right]\left[\left[t\right]\right]$. Using the actual matrix representations corresponding to these ten irreducible characters and then computing the traces, we find
\begin{equation}
\chi_1\otimes \Delta_0 ( \widehat{ H } ({\sf C}(t)) )=\frac{6}{1-4t}-\frac{6}{1-3t}+\frac{1}{1-2t},
\end{equation}
\begin{equation}
\chi_1\otimes \Delta_1 ( \widehat{ H } ({\sf C}(t)) )=1,
\end{equation}
\begin{equation}
\chi_2\otimes \Delta_0 ( \widehat{ H } ({\sf C}(t)) )= \frac{1}{1-2t},
\end{equation}
\begin{equation}
\chi_2 \otimes \Delta_1( \widehat{ H } ({\sf C}(t)) )= \frac{6}{1+2t}+1-\frac{6}{1+t},
\end{equation}
\begin{equation}
\chi_3 \otimes \Delta_0( \widehat{ H } ({\sf C}(t)) )= \frac{3}{1-3t}-\frac{1}{1-2t},
\end{equation}
\begin{equation}
\chi_3 \otimes \Delta_1( \widehat{ H } ({\sf C}(t)) )=\frac{3}{1+t}-1,
\end{equation}
\begin{equation}
\chi_4 \otimes \Delta_0( \widehat{ H } ({\sf C}(t)) )= \frac{-3}{1-t}+6,
\end{equation}
\begin{equation}
\chi_4 \otimes \Delta_1( \widehat{ H } ({\sf C}(t)) )= \frac{3}{1-t},
\end{equation}
\begin{equation}
\chi_5 \otimes \Delta_0( \widehat{ H } ({\sf C}(t)) )= \frac{3}{1-t},
\end{equation}
\begin{equation}
\chi_5\otimes \Delta_1 ( \widehat{ H } ({\sf C}(t)) )= \frac{-3}{1-t}+\frac{6}{1-2t}.
\end{equation}

Using the expression for $P$ as a linear combination of the irreducible characters, we find \begin{eqnarray}
P( \widehat{ H } ({\sf C}(t)))&=& \sum_{\stackrel{C\in \mathcal{C}_3}{C {\rm ~primitive}}} t^{ \# \textnormal{columns of }C} \nonumber \\
& =&   \frac{1}{8} \cdot  \Big(  15/(1-4t) - 18/(1-3t) +13/(1-2t) + 1 - 6/(1+t) + 3/(1+2t) \Big).\label{eq: formula}
\end{eqnarray}
The coefficient of $t^n$ of $P( \widehat{ H } ({\sf C}(t)))$ is the number of $3\times n$ Cauchon diagrams that are primitive by construction.  On the other hand, the coefficient of $t^n$ for $n\ge 1$ of the rational function 
$$ \frac{1}{8} \cdot  \Big(  15/(1-4t) - 18/(1-3t) +13/(1-2t) + 1 - 6/(1+t) + 3/(1+2t) \Big)$$ is given by $$\frac{1}{8} \cdot  \Big(  15\cdot 4^n - 18 \cdot 3^n +13 \cdot 2^n - 6\cdot(-1)^n + 3\cdot (-2)^n \Big).$$   Using equation (\ref{eq: formula}), we obtain the desired result. \qed

\section{Concluding remarks and open questions}
In this section, we make some general remarks and pose some problems we are unable to solve. We first remark that the techniques employed in \cite{BLN} in order to enumerate the $\hc$-invariant primitive ideals in $\mathcal{O}_{q}(M_{2,n})$ were more elementary than the techniques used here, and they do not extend beyond the $2\times n$ case in any obvious way.  By contrast, the techniques used here to enumerate the primitive $\hc$-invariant primes in $\mathcal{O}_q(M_{3, n})$ could in principle be used to enumerate the primitive $\hc$-invariant primes in $\mathcal{O}_q(M_{m, n})$ for any fixed $m$. The problem with larger $m$ is that the computation time grows at least exponentially in $m$ and even when $m=3$ the problem is non-trivial. We have some questions that arose during our investigations. 
\begin{quest} 

Is the image of the semigroup $\mathcal{D}_m$ under the map $$\left[\, \cdot \,\right]:\mathcal{D}_m \rightarrow \mathbb{Z}_3 \left[{\rm Ex}_3\right]/J_m$$ a group? If so what is its order?
\end{quest}

For $m=1,2,3$ this is the case and we get groups of orders $2$, $16$, and $384$ respectively.
Our next question is motivated by our calculation of the Pfaffian.  Note that if $P$ is the primitive $\hc$-invariant prime in $\mathcal{O}_q(M_{m, n})$ associated to the Cauchon diagram $C$ then Pfaffian($C$) = $\pm 2^k$. Moreover, this fact comes from looking at the zeros of an associated quadratic form $(\bmod \,2)$. Such a quantity, particularly the sign, is of great importance in knot theory---via the \emph{Arf invariant} \cite[p. 326--327]{Menasco}---and in coding theory with Reed-Muller codes \cite{quadform}. It is therefore natural to ask the following question.
\begin{quest}

The sign of the Pfaffian partitions the $\hc$-invariant primitive ideals into two classes. Is there any algebraic property that distinguishes between these two classes?
\end{quest}


\begin{thebibliography}{40}

\bibitem{reps} J. L. Alperin and R. B. Bell, {\em Groups and Representations}, Springer, New York, 1995.

\bibitem{BLN} J. Bell, S. Launois and N. Nguyen, {\em Dimension and enumeration of primitive ideals in quantum algebras}, to appear in Journal of  Algebraic Combinatorics.

\bibitem{quadform} E. R. Berlekamp, {\em Algebraic Coding Theory}, McGraw-Hill, New York, 1968.



\bibitem{bgbook} K. A. Brown and K. R. Goodearl, {\em Lectures on algebraic
  quantum groups}, Advanced Courses in Mathematics CRM Barcelona,
  Birkh\"auser, Basel, 2002.


\bibitem{Cauchon} G. Cauchon, {\em Spectre premier de $O_q \left( \mathcal{M}_n (k) \right)$, 
image canonique et s\'eparation normale}, J. Algebra {\bf 260} (2003), 519--569.

\bibitem{Cor4} S. Corteel and P. Nadeau, {\em Bijections for permutation tableaux}, to appear in the European Journal
of Combinatorics.

\bibitem{Cor1} S. Corteel and L. Williams, {\em Tableaux combinatorics for the asymmetric exclusion process}, Adv. Appl. Math. {\bf 39} (2007), 293--310.

\bibitem{Cor2} S. Corteel and L. Williams, {\em A Markov chain on permutations which projects to the asymmetric
exclusion process}, Int. Math. Res. Not. (2007), article ID mm055.



\bibitem{Cor3} S. Corteel and L. Williams, {\em A Markov chain on permutations which projects to the asymmetric
exclusion process}, posted at arXiv:0810.2916.



\bibitem{Eilenberg} S. Eilenberg, {\em Automata, Languages and Machines vol. A}, Pure and Applied Mathematics, Vol. 58. Academic Press, New York, 1974.
\bibitem{gletpams} K. R. Goodearl and E. S. Letzter, 
{\em Prime factor algebras of the coordinate ring of quantum matrices}, Proc. Amer. Math. Soc. 
121 (1994), 1017--1025.



\bibitem{GL} K. R. Goodearl and E. S. Letzter, {\em The Dixmier-Moeglin equivalence in quantum coordinate rings and quantized Weyl algebras}, Trans. Amer. Math. Soc. {\bf 352}, (2000), no. 3, 1381--1403.



\bibitem{Jo} M. Josuat-Verg\`es, {\em Bijections between pattern-avoiding fillings of Young diagrams}, posted at arXiv:0801.4928.

\bibitem{LL} S. Launois and T. H. Lenagan, {\em Primitive ideals and automorphisms of quantum matrices}, Algebr. Represent. Theory {\bf 10} (2007), no. 4, 339--365.

\bibitem{Lov} L. Lov\'asz and M. D. Plummer, {\em Matching theory}, Ann. Discrete Math. {\bf 29}, North-Holland, 1986.

\bibitem{Meduna} A. Meduna, {\em Automata and Languages Theory and Applications}, Springer-Verlag, Ltd., London, 2000.

\bibitem{Menasco} W. Menasco and M. Thistlethwaite, {\em Handbook of knot theory}, Elsevier B. V., Amsterdam, 2005.

\bibitem{Post} A. Postnikov, {\em Total positivity, Grassmannians, and networks}, posted at arXiv:math/0609764.

\bibitem{Sipser} M. Sipser, {\em Introduction to the Theory of Computation}, Second Edition, Thomson Course Technology, Boston, 2006.

\bibitem{Ste} E. Steingrimsson and L. Williams, {\em Permutation tableaux and permutation patterns}, J. Comb. Th. A {\bf 114} (2007), 211-234.

\bibitem{Wi} L. Williams, {\em Enumeration of totally positive Grassmann cells}, Adv. Math. {\bf 190} (2005), 319--342.

\end{thebibliography}
\end{document}